%% file: journal_arxiv.tex
\documentclass[draftcls,onecolumn]{IEEEtran}

\makeatletter
\let\proof\@undefined
\let\endproof\@undefined
\makeatother

\usepackage{amsmath}
\usepackage{amsfonts}
\usepackage{amssymb}
\usepackage{mathrsfs}
\usepackage{fix2col}
\usepackage{latexsym}
\usepackage[dvips]{epsfig}

\usepackage{cite}
\usepackage{amsmath, amsthm, amssymb,graphicx,bbm}
\usepackage{verbatim}
\usepackage{psfrag}
\usepackage{algorithm2e}
\usepackage{bm}
\usepackage{color}
\usepackage[tight,footnotesize]{subfigure}

\input{header.tex}
\newcommand{\ch}{{}} 
\begin{document}
\title{\Large \bf Geometry of Power Flows and Optimization in Distribution Networks }
 \author{Javad Lavaei, David Tse and Baosen Zhang
 \thanks{Authors sorted alphabetically, all three contributed equally to this work.}%
\thanks{Javad Lavaei is  with  the Department of Electrical Engineering, Columbia University (email: lavaei@ee.columbia.edu). David Tse is with the Department of Electrical Engineering and Computer Sciences, University of California, Berkeley  (email:dtse@eecs.berkeley.edu). Baosen Zhang is the Department of Civil and Environmental Engineering and Management Science \& Engineering at Stanford University (email:baosen.zhang@gmail.com). }%
\thanks{The work of B. Zhang and D. Tse was supported in part by the National Science Foundation (NSF) under grant CCF-0830796. B. Zhang was also supported by a National Sciences and Engineering Research Council of Canada  Postgraduate scholarship.}}
\maketitle
\begin{abstract}
We investigate the geometry of injection regions and its relationship to optimization of power flows in tree networks. The injection region is the set of all vectors of bus power injections that satisfy the network and operation constraints. The geometrical object of interest is the set of Pareto-optimal points of the injection region.  If the voltage magnitudes are fixed, the injection region of a tree network can be written as a linear transformation of the product of two-bus injection regions, one for each line in the network. Using this decomposition, we show that under the practical condition that the angle difference across each line is not too large, the set of Pareto-optimal points of the injection region remains unchanged by taking the convex hull. { Moreover, the resulting convexified  optimal power flow problem can be efficiently solved via }{ semi-definite programming or second order cone relaxations. These results improve upon earlier works by removing the assumptions on active power lower bounds.} It is also shown that our practical angle assumption guarantees two other properties: (i) the uniqueness of the solution of the power flow problem, and (ii) the non-negativity of the locational marginal prices. Partial results are presented for the case when the voltage magnitudes are not fixed but can lie within certain bounds.

\end{abstract}
\section{Introduction}


AC optimal power flow (OPF) is a basic problem in power engineering. The problem is to efficiently allocate power in the electrical network, under various operation constraints on voltages, flows, thermal dissipation and bus powers. For general networks, the OPF problem is known to be nonconvex  and is challenging \cite{Huneault91,Ian_1}. Some of the earlier analysis has focused on understanding the existence and the behavior of load flow around local solutions \cite{load_flow_1,load_flow_3}. Recently, different convex relaxation techniques have been applied to the OPF problem in an attempt to find global solutions \cite{opt_j_4,opf_sdp1}. It was recently observed in\cite{javad_zero1}  that many practical instances of the OPF problem can be convexified via a rank relaxation. This observation spurs the question:  when can an OPF problem be convexifed and solved efficiently?
This question was partially answered in several recent independent works \cite{Zhang12,Lavaei11c,Bose11}: the convexification of OPF is possible if the network has a tree topology and  some conditions on the bus power constraints hold.  The goal of this paper is to provide a unified understanding of these results through a deeper investigation of the underlying geometry of the optimization problem. Through this understanding, we are also able to strengthen these earlier results.

There are three reasons why it is worthwhile to focus on tree networks. First, although OPF is traditionally solved for transmission networks, there is an increasing interest in optimizing power flows in distribution networks due to the emergence of demand response and distributed generation \cite{DOE}.  Unlike transmission networks, most distribution networks have a tree topology. {\ch Second, as will become apparent, assuming a tree topology is a natural simplification of the general OPF problem (\cite{load_flow_2} made a similar observation for the power flow problem), and results for this simplified problem will shed light on the general problem.}  Third, as was shown in \cite{Lavaei11c}, if one is allowed to put phase shifters in the network, then the OPF problem for general network topologies can essentially be reduced to one for tree networks.  


{ Following \cite{Zhang12}, our approach to the problem is based on an investigation of the convexity properties of the {\em power injection region}. The injection region is  the set of all vectors of feasible real power injections $P_i$'s (both generations and withdraws) at the various buses that satisfy the given network and operation constraints. We are particularly interested in the {\em Pareto-front} of the injection region; these are the points on the boundary of the region for which one cannot decrease any component without increasing another component. The significance of the Pareto-front  is that the optimal solution of OPF problems with increasing objective functions defined on the injection region must lie there.
The first question we are after is: although the injection region is nonconvex, when does its Pareto-front remain unchanged upon taking the convex hull of the injection region? This property would ensure that any OPF problem over the injection region is  {\em convexifiable}: solving it over the larger convex hull of the injection region would yield an optimal solution to the original nonconvex problem  (in the sequel, we will abbreviate this property by simply saying that the Pareto-front is convex).  While convexifiability is a desirable intrinsic property of any optimization problem, a second question of interest is: can the resulting convexified OPF problem be solved efficiently?
}

\begin{figure}[ht]%
\centering
\includegraphics[width=2.5cm,height=2cm]{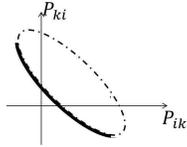}
\caption{The feasible set for the two flows along a line when there are power flow constraints. It is a subset of an ellipse which is the feasible set when there are no constraints other than the fixed voltage magnitudes at the two buses. In this example, the feasible set is part of the Pareto-front of the ellipse.}
\label{fig:flow_region}%
\end{figure}

To answer these questions, the first step is to view the injection region as a linear transformation of the higher dimensional {\em power flow} region: this is the set of all vectors of feasible real power flows $P_{ik}$'s, one along each direction of each line. We first focus on the case when the voltage magnitudes are fixed at all buses. In the space of power flows, the network and operation constraints in general networks can be grouped into three types:
\begin{enumerate}
\item {\em local} constraints on the two flows $P_{ik}$ and $P_{ki}$ along each line $(i,k)$: these include angle, line flow and thermal constraints. Figure \ref{fig:flow_region} gives an example of the (nonlinear) feasible set of $(P_{ik},P_{ki})$ due to flow constraints.  Note that all these local constraints are effectively angle constraints. Also this feasible set can be interpreted as the power flow region of a two-bus network with impedance given by that of the line $(i,k)$.
\item {\em global} constraints on the flows due to bus power constraints. These are linear constraints.
\item {\em global} Kirchhoff constraints on the flows due to cycles. These are non-linear constraints.
\end{enumerate}

The third type of constraints is most complex; they are global and non-linear. By focusing on tree networks, we are left only with constraints of type $1$ and type $2$.  In this case, it is easy to see that the power flows along different lines are {\em decoupled}, save for the global (but linear) bus power constraints. The overall power flow region is thus simply the {\em product} of the two-bus power flow regions, one for each line, intersecting the bus power constraints.

{ By exploiting this geometric structure of the power flow region, we answer the first question posed above: if the two-bus power flow region associated with each line  itself has a convex Pareto-front, then  the overall injection region has a convex Pareto-front.  Thus, a {\em local} convexity property guarantees a {\em global} convexity property. It is shown that the local property holds whenever the angle difference along every line is constrained to be not too large (say less than $45^\circ$). Note that the angle constraints are not an additional constraints in the OPF problem, because as we will argue, existing constraints due to line flow or thermal constraints can be thought as angle constraints. 

Concerning the second computational question,  we observe that in our geometric picture,  a semi-definite programming relaxation corresponds to taking the convex hull of the purely voltage-constrained injection region followed by intersection with the local and global constraints.  This convex relaxation results in general in a set larger than the convex hull of the injection region. However, it turns out that our analysis of the first question in fact implies that the Pareto-front remains unchanged even with this more relaxed convexification. This provides a resolution to the computational question.
}

The present work improves upon the earlier papers on tree networks  \cite{Zhang12,Lavaei11c,Bose11,SGC2011,CDC2012} in two ways. First, the arguments used in \cite{Zhang12,Lavaei11c,Bose11} are  algebraic and some used non-trivial matrix fitting results and \cite{SGC2011,CDC2012} use algebraic SOCP relaxations, while the present paper uses entirely elementary geometric arguments. This geometric approach provides much more insight on the roles of the various types of constraints and also explains how the assumption of tree topology simplifies the problem. Second, the convexity results in all of the earlier papers require some restriction on the bus power lower bounds (no lower bounds are allowed in \cite{Lavaei11c,Bose11}, and any two buses that are connected cannot simultaneously have lower bounds in \cite{Zhang12}). The results in this paper require no such conditions. Instead, they are replaced by constraints on the differences between voltage angles at adjacent buses,  which we verify to be satisfied in practice. This latter condition imposes a local constraint on the power flows, and is discovered through the geometric decomposition of the power flow region. We note that not all OPF problem on trees have objective functions on the power flow region. For example, the conservation voltage reduction problem is to reduce the voltage of the buses to some threshold \cite{CDC2012}, and the feasible space of voltages and is an interesting problem of future study.

We also show that the angle assumption gives rise to two other important properties: (i) the solution of a power flow problem becomes unique, and  (ii) the locational marginal prices (LMPs) are never negative when all bid functions are positive. The LMP is a common pricing signal used in practice for charging customers and paying generators located at various buses. It is known that in congested transmission networks LMP at some buses could be negative \cite{kirschen04} even when all bid functions are positive. We show that this situation does not happen for a tree network with realistic angle assumptions.

{ The paper is organized as follows. In Section \ref{sec:notation}, we state the physical model used in the paper. Section \ref{sec:O1} focuses on the case when the voltages magnitudes at all buses are fixed. We first start with a two-bus network with angle, thermal and flow constraints. Then, we consider a general tree network with only local constraints. Finally, we add the global bus power constraints to arrive at our full results. We also study the implications of our result on the uniqueness of power flow solutions and the non-negativity of locational marginal prices. Section \ref{sec:variable} extends some of the results to networks with variable voltage magnitudes. Section \ref{sec:simu} shows simulation results that validate our theoretical insights, and Section \ref{sec:con} concludes the paper. Similar results can be derived if the network has reactive power constraints in addition to active power constraints. Due to space limitations, we do not include them here.
}

Summary of notations we use throughout the paper:
\begin{itemize} 
\item {\bf vectors and matrices:}   We use the notations $\x$ and $\bd{X}$ to denote vectors and matrices, respectively.
Given two real vectors $\bd{x}$ and
$\bd{y}$ of the same dimension, the notation $\bd{x} \leq \bd{y}$
denotes a component-wise inequality. We denote Hermitian transpose of a matrix with $(\cdot)^H$ and conjugation with $\conj(\cdot)$. The vector $\x \odot \y$ is the component-wise product of the two vectors $\x$ and $\y$, and $\diag(\bd{X})$ returns the vector containing the diagonal elements of the matrix $\bd{X}$. The notation $j$ is reserved for $\sqrt{-1}$ in this work.
\item {\bf sets:} We use scripted capital letters $\mc{A},\mc{B},\dots$ to represents sets, which are assumed to be subsets of $\R^{n}$ unless otherwise stated. Given a set $\mathcal{A}$, $\conv(\mathcal{A})$ denotes the convex hull of $\mathcal{A}$. A point $\x \in \mc{A}$ is Pareto-optimal if there does not exist another point $\y \in \mc{A}$ such that $\y \leq \x$ with strict inequality in at least one coordinate. Let $\mc{O}(\mc{A})$ denote the set of all Pareto-optimal points of $\mc{A}$, which is sometimes called the Pareto front of $\mc{A}$. Note that if a strictly increasing function is minimized over $\mc{A}$, its optimal solution must belong to $\mc{O}(\mc{A})$. 
\end{itemize}

\section{Model } \label{sec:notation}
Consider an AC electrical power network with $n$ buses. With no loss of generality, we assume that the network is a connected graph.  Following the convention in power engineering, complex scalars representing voltage, current and power are denoted by capital letters.  We write $i \sim k$ if bus $i$ is connected to $k$, and $i \nsim k$ if they are not connected. We often regard the network as a graph with the vertex set $\mathcal V=\{1,\dots,n\}$ and the edge set $\mc{E}$. For example, the notation $(i,k) \in \mc{E}$ implies that there exists a line connecting bus $i$ and bus $k$. Let $z_{ik}$ denote the complex impedance of the line $(i,k)$  and $y_{ik}=\frac{1}{z_{ik}}=g_{ik}-jb_{ik}$ represent its admittance, where $g_{ik},b_{ik}\geq0$. Define the admittance matrix  $\bd{Y}$ as
\begin{equation} \label{eqn:Y}
Y_{ik}=
\begin{cases}
\sum_{l \sim i} y_{il} +y_{ii} & \mbox{ if } i=k \\
-y_{ik} & \mbox { if } i \sim k \\
0 & \mbox { if } i \nsim k
\end{cases},
\end{equation}
where $y_{ii}$ is the shunt admittance to ground at bus $i$. Note that this matrix is symmetric.

Let $\bd{v}=(V_1,V_2,\dots,V_n) \in \C^n$ be the vector of complex bus voltages and $\bd{i}=(I_1,I_2,\dots,I_n) \in \C^n$ be the vector of complex currents, where $I_i$ is the total current flowing out of bus $i$ to the rest of the network. By Ohm's law and Kirchoff's Current Law, $\bd{i}=\bd{Y}\bd{v}$.  The complex power injected to bus $i$ is equal to $S_i=P_i+j Q_i=V_i I_i^H$, where $P_i$  and $Q_i$ denote the net active and reactive powers at this bus, respectively. Let $\bd{p}=(P_1,P_2,\dots,P_n)$ be the vector of real powers, which can be written as $\bd{p}=\Real(\bd{v} \odot \conj(\bd{i}))=\Real(\bd{v} \odot (\bd{Y} \bd{v}))=\Real(\diag(\bd{v} \bd{v}^H \bd{Y}^H))$.

\section{Fixed Voltage Pareto Optimal Points} \label{sec:O1}
\subsection{Two-Bus Network With Angle, Thermal and Flow Constraints}
\begin{figure}[ht]
\centering
\psfrag{V1}{$V_1$}
\psfrag{V2}{$V_2$}
\psfrag{P1}{$P_1$}
\psfrag{P2}{$P_2$}
\psfrag{g-bj}{$g-j b$}
\includegraphics[scale=0.8]{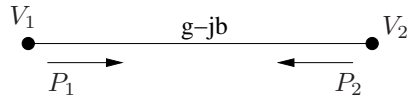}%
\caption{A two-bus network.}%
\label{fig:2bus}%
\end{figure}

\begin{figure}[ht]%
\centering
\psfrag{P1}{$P_1$}
\psfrag{P2}{$P_2$}\hspace{-5mm}
\subfigure[lossy]{
\includegraphics[scale=0.1]{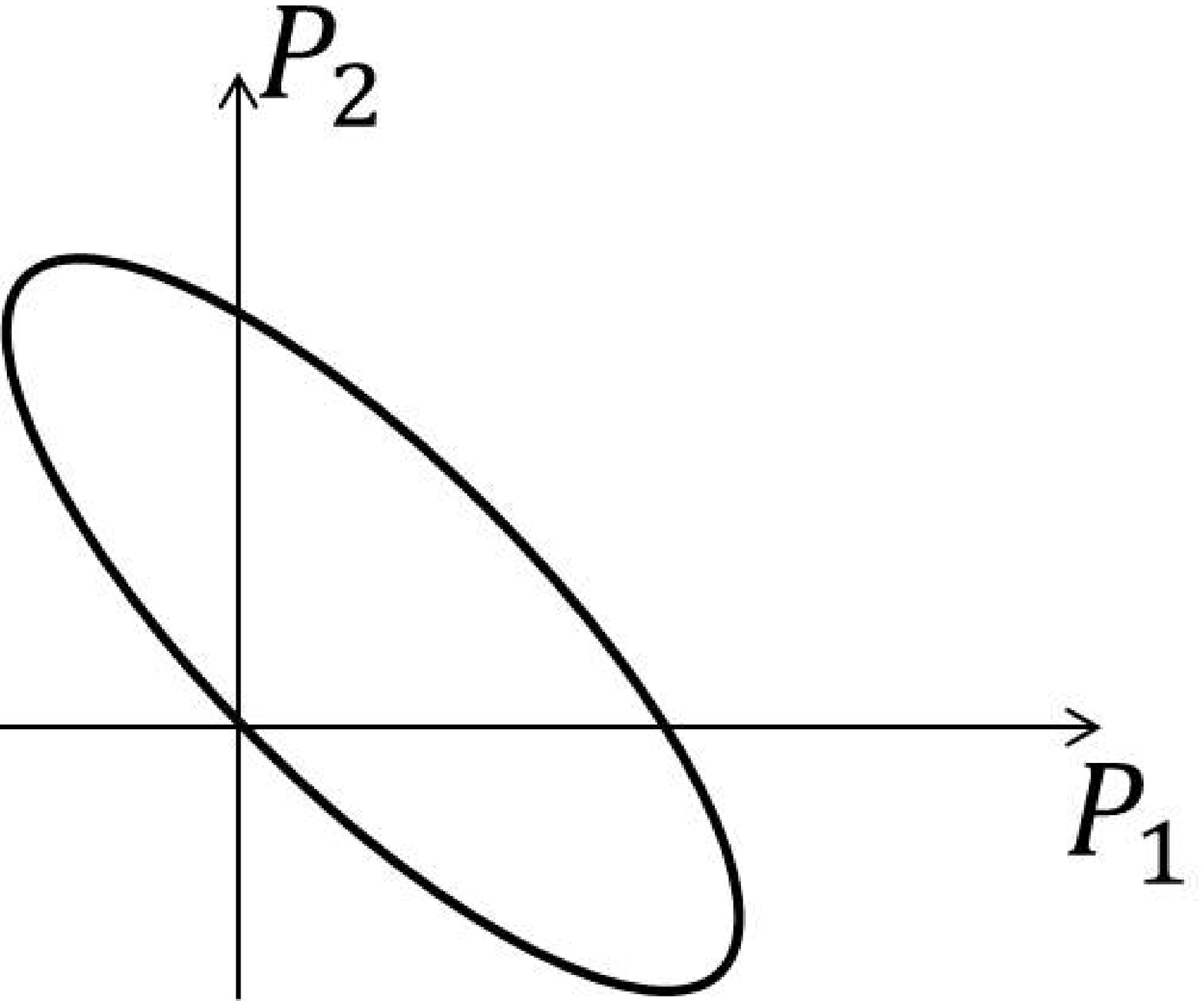} \label{fig:P12b}}\hspace{5mm}
\subfigure[lossless]{
\includegraphics[scale=0.1]{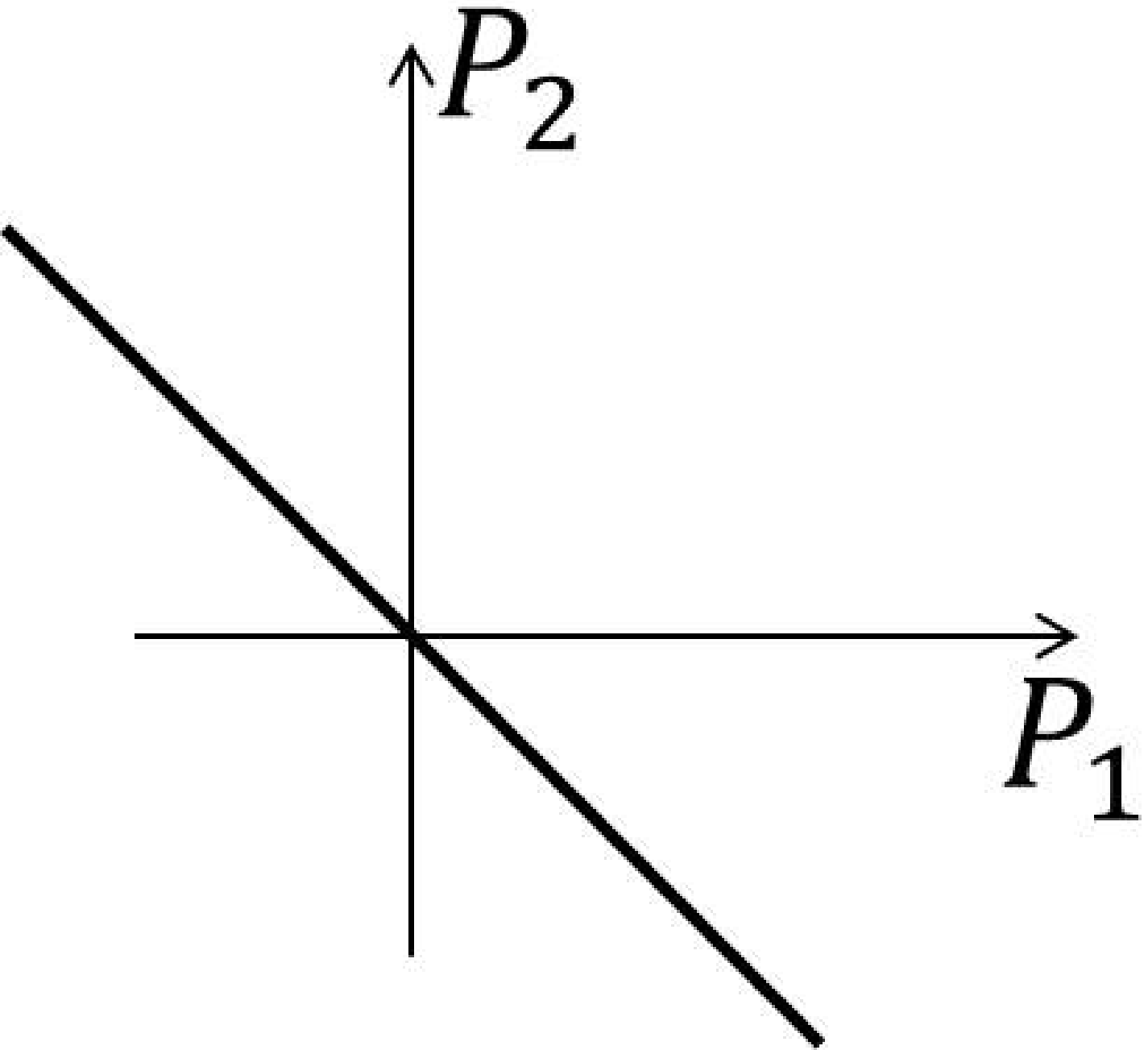} \label{fig:P12g}}
\caption{The region  defined by \eqref{eqn:P}: (a) shows the region corresponding to $|V_1|=|V_2|=1$ (per unit), $b=5$ and $g=1$; (b) shows the region for a lossless line. }%
\label{fig:ellipse}%
\vspace{-4mm}
\end{figure}

Consider the two-bus network in Figure \ref{fig:2bus} with the line admittance $g-j b$.
Let the complex voltages at buses 1 and 2 be expressed as $V_1=|V_1| \exp (j \te_1)$ and $V_2=|V_2| \exp (j \te_2)$. Throughout this subsection, assume that the magnitudes $|V_1|$ and $|V_2|$ are fixed, while $\theta_1$ and $\theta_2$ are variable. The power injections at the two buses are given by
\begin{subequations} \label{eqn:P}
\begin{align}
P_{1} &=|V_1|^2 g+ |V_1||V_2| b \sin (\te) -|V_1||V_2| g \cos(\te) \\
P_{2} &=|V_2|^2 g- |V_1||V_2| b \sin (\te) -|V_1||V_2| g \cos(\te),
\end{align}
\end{subequations}
where $\te=\te_1-\te_2$.
Since the network has only two buses, $P_1=P_{12}$ and $P_2=P_{21}$, where $P_{ik}$ is the power flowing out of bus $i$ to bus $k$. Since the voltage magnitudes are fixed, the  power flows between the buses can both be described in terms of the single parameter $\te$. Notice that a circle centered at the origin and of radius $1$ can be parameterized as $(\cos(\te),\sin(\te))$. Therefore, \eqref{eqn:P} represents an affine transformation of a circle, which leads to an ellipse. This ellipse contains all points $(P_{1},P_{2})$ satisfying the inequality
\begin{equation}
\nonumber
\bigg\|\left[\begin{array}{cc} b&-g\\ -b& -g\end{array}\right]^{-1}\left[\begin{array}{cc}P_{1}-|V_1|^2g\\P_{2}-|V_2|^2g\end{array}\right]\bigg\|_2 = |V_1||V_2|,
\end{equation}
where $\|\cdot\|_2$ denotes the 2-norm operator. As can be seen from the above relation, the ellipse is centered at $(|V_1|^2 g, |V_2|^2 g)$, where its major axis is at an angle of $-45^{\circ}$ to the x-axis with length $|V_1V_2| b$ and its minor principal axis has length $|V_1 V_2| g$. {\ch If the line is lossy, the injection region is a hollow ellipse as shown in Figure \ref{fig:P12b}.  If the line is lossless, the ellipse is degenerate and collapses into a line through the origin as shown in Figure \ref{fig:P12g}.} 
{\ch In practice most lines in distribution networks are lossy with $b/g$ ratio typically between 0.2 to 5 (instead of $>10$ in transmission networks), with overhead cable having higher $b/g$ ratios.  \cite{testfeeders,Kersting06}. Thus, the interesting and practical case is when the region is a hollow ellipse. Note that the convex hull of this region is the filled ellipse.}

Now, we investigate the effect of thermal, line flow and angle constraints. Since the network has fixed voltage magnitudes, the thermal loss and line flow constraints can be recast as angle constraints of the form $\ul{\te} \leq \te \leq \ov{\te}$ for some limits $\ul{\te}\in[-\pi,0]$ and $\ov{\te} \in[0,\pi]$. More precisely, the loss of the line, denoted by $L_{12}$, can be calculated as
\begin{equation} \label{eqn:loss}
\begin{aligned}
L_{12}&= |V_1 -V_2|^2 g=P_{12}+P_{21}\\
&=|V_1|^2g-2|V_1||V_2|g \sin (\te)+ |V_2|^2 g.
\end{aligned}
\end{equation}
It follows from the above equality that a loss constraint  $L_{12}\leq \ov{L}_{12}$ (for a given $\ov{L}_{12}$) can be translated into an angle constraint.  Likewise, the line flow inequalities $P_{12} \leq \ov{P}_{12}$ and $P_{21} \leq \ov{P}_{21}$ are also angle constraints. As a result, we restrict our attention only to angle constraints in the rest of this part.

We define the injection region $\mathcal P$ to be the set of all points $\{(P_1,P_2)\}$ given by \eqref{eqn:P} by varying $\te \in [ \ul{\te},\ov{\te}]$.  The bold curve in Figure \ref{fig:flow_region}
represents the injection region after a certain angle constraint.

The key property of the non-convex feasible set $\mc P$ for a two-bus network is that Pareto front of P is the same as the Pareto-front of the convex hull of $\mc P$ (see  Figure \ref{fig:flow_region}). To understand the usefulness of this property in solving an optimization problem over this region, consider the following pair of optimization problems for a strictly increasing function $f$:
\begin{subequations} \label{eqn:opf2}
\begin{align}
\mbox{minimize } & f(P_1,P_2) \\
\mbox{subject to } & (P_1,P_2) \in \mc{P},
\end{align}
\end{subequations}
and
\begin{subequations} \label{eqn:opf2_conv}
\begin{align}
\mbox{minimize } & f(P_1,P_2) \\
\mbox{subject to } & (P_1,P_2) \in \conv(\mc{P}).
\end{align}
\end{subequations}
Since $f$ is strictly increasing in both of its arguments, the optimal solution to \eqref{eqn:opf2_conv} must be on the Pareto boundary of the feasible set; therefore both optimization problems share the same solution $(P_1^*,P_2^*) \in \mc{P}$. This implies that instead of solving the non-convex problem \eqref{eqn:opf2}, one can equivalently solve the  optimization (\ref{eqn:opf2_conv}) that is always convex for a convex function $f$. Hence, even though  $\mc{P}$ is not convex, optimization over $\mc{P}$ and $\conv(\mc{P})$ is equivalent for a broad range of optimization problems due to the following lemma.
\begin{lem} \label{lem:O_nocon}
Let $\mc{P}\in \R^2$ be the two-bus injection region defined in \eqref{eqn:P} by varying $\theta$ over $[\ul{\te},\ov{\te}]$. The relation $\mc{O}(\mc{P})=\mc{O}(\conv(\mc{P}))$ holds.
\end{lem}

\subsection{General Network With Local Constraints}
In this subsection, we extend Lemma \ref{lem:O_nocon} to an arbitrary tree network with local constraints while section \ref{sec:O} states and proves the general result with both local and global constraints.

First, we express the injection region of a general tree as a \textit{linear transformation} of the power flow region.
Given a general network described by its admittance matrix $\bd{Y}$, consider a connected pair of buses $i$ and $k$. Let $P_{ik}$ denote the power flowing from bus $i$ to bus $k$ through the line $(i,k)$ and $P_{ki}$ denote the power flowing from bus $k$ to bus $i$. Similar to the two-bus case studied earlier, one can write:
\begin{subequations}
\nonumber
\begin{align}
P_{ik} &=|V_i|^2 g_{ik}+ |V_i||V_k| b_{ik} \sin \te_{ik} -|V_i||V_k| g_{ik} \cos \te_{ik} \nn \\
P_{ki} &=|V_k|^2 g_{ik}- |V_i||V_k| b_{ik} \sin \te_{ik} -|V_i||V_k| g_{ik} \cos \te_{ik}, \nn
\end{align}
\end{subequations}
where $\te_{ik}=\te_i-\te_k$. The tuple $(P_{ik},P_{ki})$ is referred to as the flow on the line $(i,k)$. As in the two-bus case, all the thermal and line flow constraints can be cast as a constraint on the angle $\te_{ik}$. Note that the angle constraint on $\te_{ik}$ only affects the flow on the line $(i,k)$; therefore it is called a local constraint.

There are $2|\mathcal E|$ numbers describing the flows in the network. Let $\mc{F}$ denote the feasible set of the flows in $\R^{2|\mathcal E|}$, where the bus voltage magnitudes are fixed across the network and each flow satisfies its local constraints. Recall that the net injection at bus $i$ is related to the line flows through the relation $P_i = \sum_{k:k \sim i} P_{ik}$.
This motivates the introduction of an $n \times 2|\mathcal E|$ matrix $\bd{A}$ defined below with rows indexed by the buses and the columns indexed by the lines:
\begin{equation} \label{eqn:A}
A(i,(k,l))=
\begin{cases}
1 & \mbox{ if } i=k \\
0 & \mbox{ otherwise.}
\end{cases}
\end{equation}
The matrix $\bd{A}$ can be seen as a generalization of the edge-to-node adjacency matrix of the graph. The injection vector $\bd{p} \in \R^n$ and the flow vector $\bd{f} \in \mc{F}$ are related by $\p=\bd{A} \bd{f}$. We express the set of line flows as $\{P_{ik},P_{ki}\}$ and say that $\p$ is achieved by the set of flows.
This implies that the feasible injection region $\mc{P}$ is given by
\begin{equation} \label{eqn:AF}
\mc{P} = \bd{A} \mc{F}.
\end{equation}
Since the above mapping is linear, it is straightforward to show that $\conv(\mc{P})= \bd{A} \conv(\mc{F})$.

We now demonstrate that $\mc{F}$  has a very simple structure: it is simply a \textit{product} of the two-bus flow regions, one for each line in the network:
\begin{equation} \label{eqn:F}
\mc{F}=\prod_{(i,k) \in \mc{E}} \mc{F}_{ik},
\end{equation}
 where the two-dimensional set $\mc{F}_{ik}$ is the two-bus flow region of the line $(i,k)$.
In other words, the flows along different lines are decoupled. To substantiate this fact, it suffices to show that the  flow on an arbitrary line of the network can be adjusted without affecting the flows on other lines. To this end, consider the  line $(i,k)$ and a set of voltages with the angles $\te_1,\dots,\te_n$.  The power flow along the line $(i,k)$ is a function of $\te_{ik}=\te_i-\te_k$. Assume that we want to achieve a new flow on the line associated with some angle $\tilde{\te}_{ik}$. In light of the tree structure of the network, it is possible to find a new set of angles $\tilde{\te}_1,\dots,\tilde{\te}_n$ such that $\tilde{\te}_i-\tilde{\te}_k=\tilde{\te}_{ik}$ and that the angle difference is preserved for every line in $\mathcal E\backslash (i,k)$.

Due to this product structure of $\mc{F}$, it is possible to generalize  Lemma \ref{lem:O_nocon}.
\begin{lem} \label{lem:O_tree}
Given a tree network with fixed voltage magnitudes and local angle constraints, consider the  injection set $\mc{P}$  defined  in \eqref{eqn:AF}. The relation $\mc{O}(\mc{P})=\mc{O}(\conv(\mc{P}))$ holds.
\end{lem}

\emph{Proof:}
First, we show that $\mc{O}(\conv(\mc{P}))\subseteq \mc{O}(\mc{P})$. Given $\p \in \mc{O}(\conv(\mc{P}))$, let $\{(P_{ik},P_{ki})\} \in \conv(\mc{F})$ be the set of flows that achieves $\p$. Consider a line $(i,k)\in\mathcal E$. Since $\mc{F}$ is a product space and $\p \in \mc{O}(\conv(\mc{P}))$, we have $(P_{ik},P_{ki}) \in \mc{O}(\conv(\mc{F}_{ik}))$. Moreover, it follows from  Lemma \ref{lem:O_nocon} that $\mc{O}(\conv(\mc{F}_{ik}))=\mc{O}(\mc{F}_{ik})$. Therefore, $(P_{ik},P_{ki}) \in \mc{F}_{ik}$ for every line $(i,k)$. This gives $\p \in \mc{P}$ and consequently $\p \in \mc{O}(\mc{P})$.

Next, we show that $\mc{O}(\mc{P}) \subseteq \mc{O}(\conv(\mc{P}))$. Given $\p \in \mc{O}(\mc{P})$,  assume  that $\p \notin \mc{O}(\conv(\mc{P}))$. Then, there exists a point $\p' \in \mc{O}(\conv(\mc{P}))$ such that $\p' \leq \p$ with strict inequality in at least one coordinate. By the first part of the proof, we have $\p' \in \mc{P}$, which contradicts $\p \in \mc{O}(\mc{P})$.\hfill$\blacksquare$

\subsection{Two-Bus Network with Bus Constraints}
\begin{figure}
\centering
\subfigure[]{
\includegraphics[width=2cm]{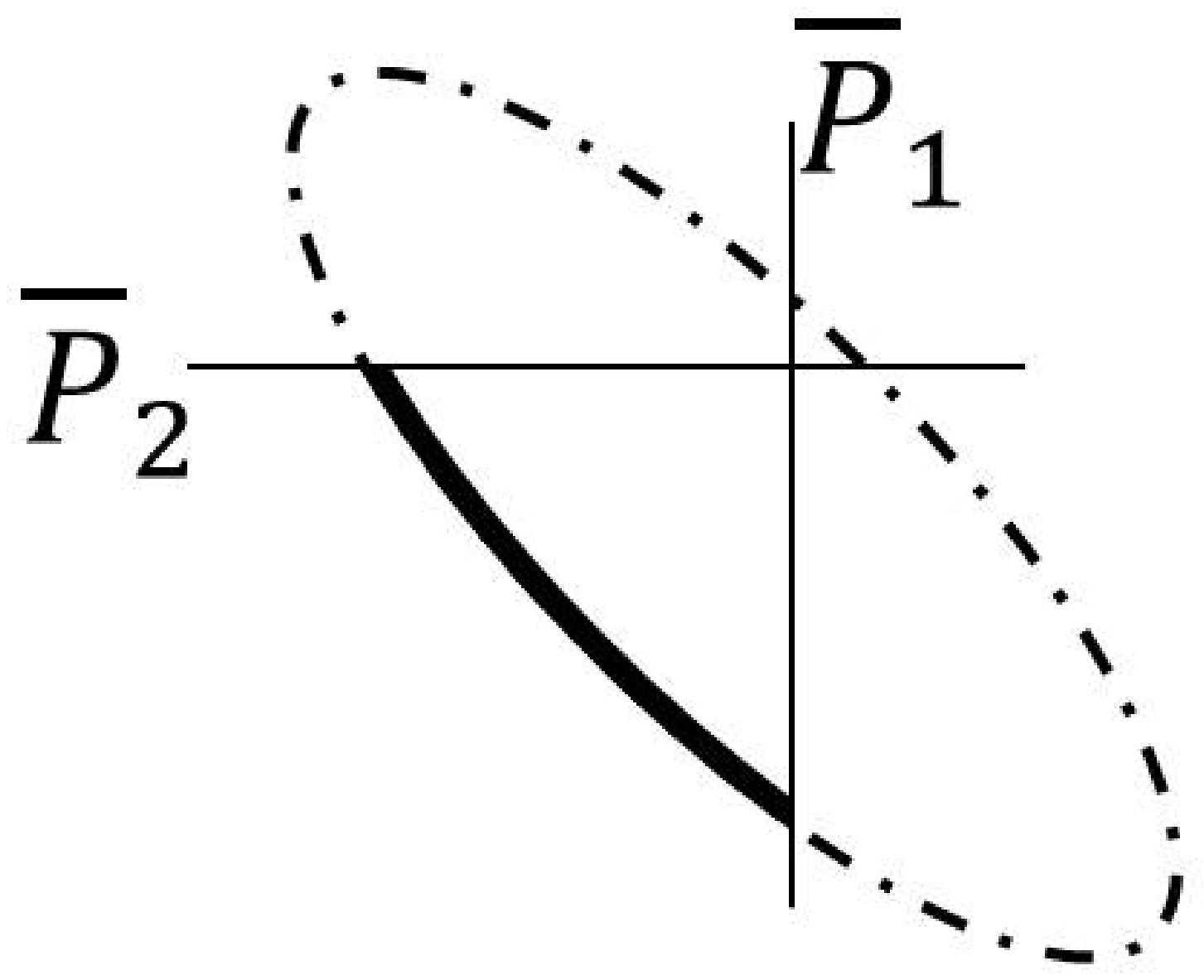}%
\label{fig:2busP1}}
\subfigure[]{
\includegraphics[width=2cm]{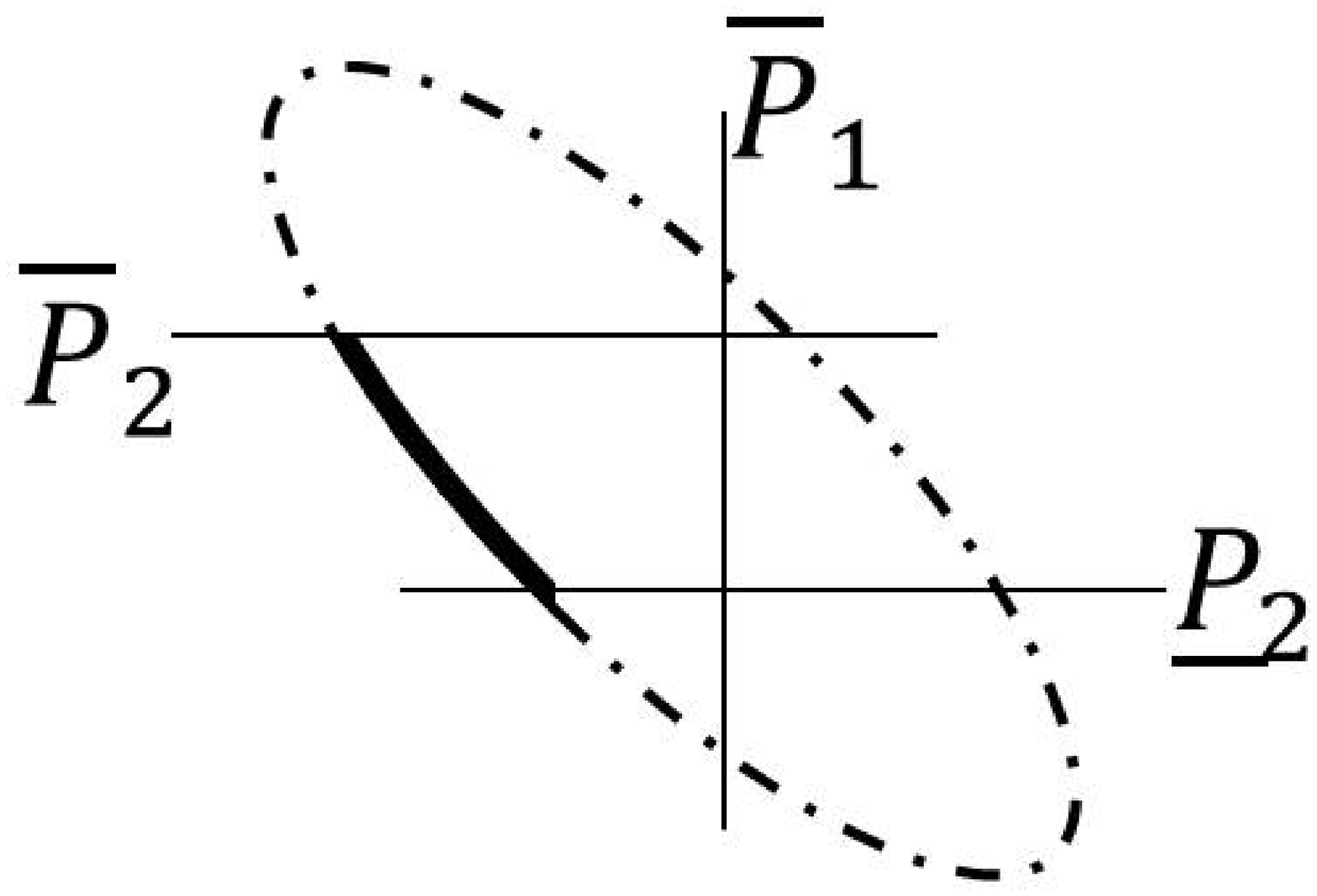}
\label{fig:2busP2}}
\subfigure[]{
\includegraphics[width=2cm]{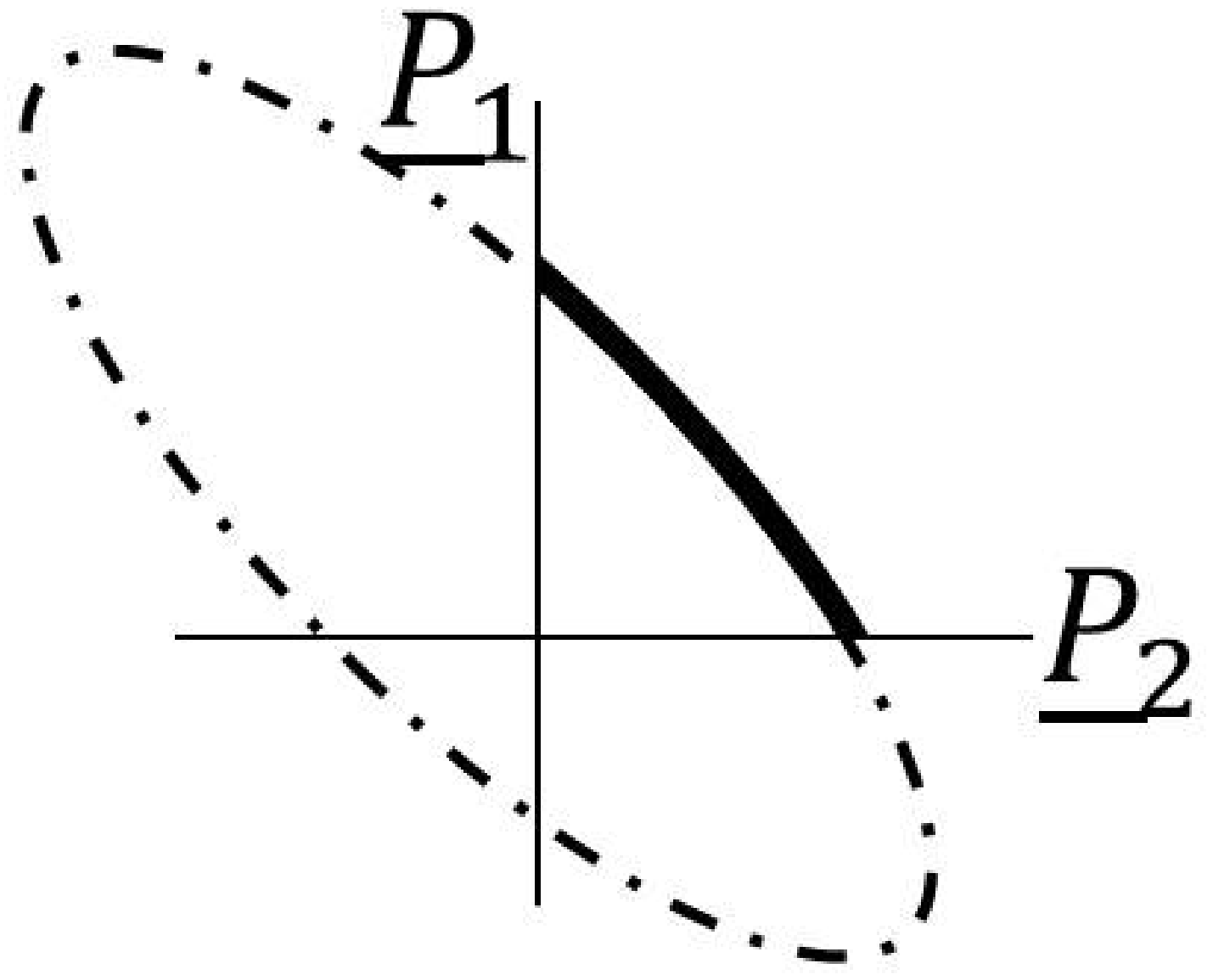}
\label{fig:2busP3}}
\caption{Three possible cases for the bus power constrained injection region. }
\label{fig:2busP}
\end{figure}

So far, we have studied tree networks with local angle constraints and without global bus power constraints. We want to investigate the effect of bus power constraints. We first consider the two-bus network shown in Figure \ref{fig:2bus}, and incorporate the angle constraints together with the bus active power constraints of the form $\ul{P}_i\leq P_i \leq \ov{P}_i$ for $i=1,2$. Let $\mc{P}_\te=\{(P_1,P_2): |V_1|=\ov{V}_1,\ |V_2|=\ov{V}_2,\ \ul{\te}\leq \te \leq \ov{\te}\}$ be the angle-constrained injection region, and $\mc{P}_P=\{(P_1,P_2):\ul{P}_i \leq P_i \leq \ov{P}_i,\; i=1,2\}$ be the bus power constrained region, where $\ov{V}_1$ and $\ov{V}_2$ are the given nominal values of the voltage magnitudes. The overall injection region is given by the intersection of the two regions through the equation
\begin{equation}
\mc{P}=\mc{P}_\te \cap \mc{P}_P.
\end{equation}
 There are several possibilities for the shape of $\mc{P}$, as visualized in Figures
\ref{fig:2busP1}, \ref{fig:2busP2} and \ref{fig:2busP3}. In Figure
\ref{fig:2busP1}, both buses have power upper bounds. In Figure \ref{fig:2busP2}, $P_1$ has upper bound, while
$P_2$ has both upper and lower bounds. In Figure \ref{fig:2busP3}, both buses have lower bounds. It can be observed that $\mc{O}(\mc{\mathcal P})=\mc{O}(\conv(\mc{\mathcal P}))$ for Figures \ref{fig:2busP1}-\ref{fig:2busP2}, but this desirable property does not hold for Figure \ref{fig:2busP3}.

Figure \ref{fig:2busP3} means that in the presence of active power lower bounds, the relationship $\mc{O} (\mc{P}) = \mc{O}(\conv(\mc{P}))$ does not always hold. This is the reason for the various assumptions made about bus power lower bounds in \cite{Zhang12,Lavaei11c,Bose11}. Note that $\te_{12}$ in Figure \ref{fig:2busP3} is allowed to vary from $-\pi$ to $\pi$. However, the angles are often constrained in practice by thermal and/or stability conditions. For example, the thermal constraints usually limit the angle difference on a line to be less than $10^{\circ}$. Figure \ref{fig:thermal} in the appendix shows a typical distribution network together with its thermal constraints, from which it can be observed that each angle difference is restricted to be less than $7^{\circ}$. Flow constraints also limit the angle differences in a similar fashion.

\begin{figure}%
\centering
\subfigure[Feasible]{ \label{fig:2bus_theta1}
\includegraphics[width=2cm]{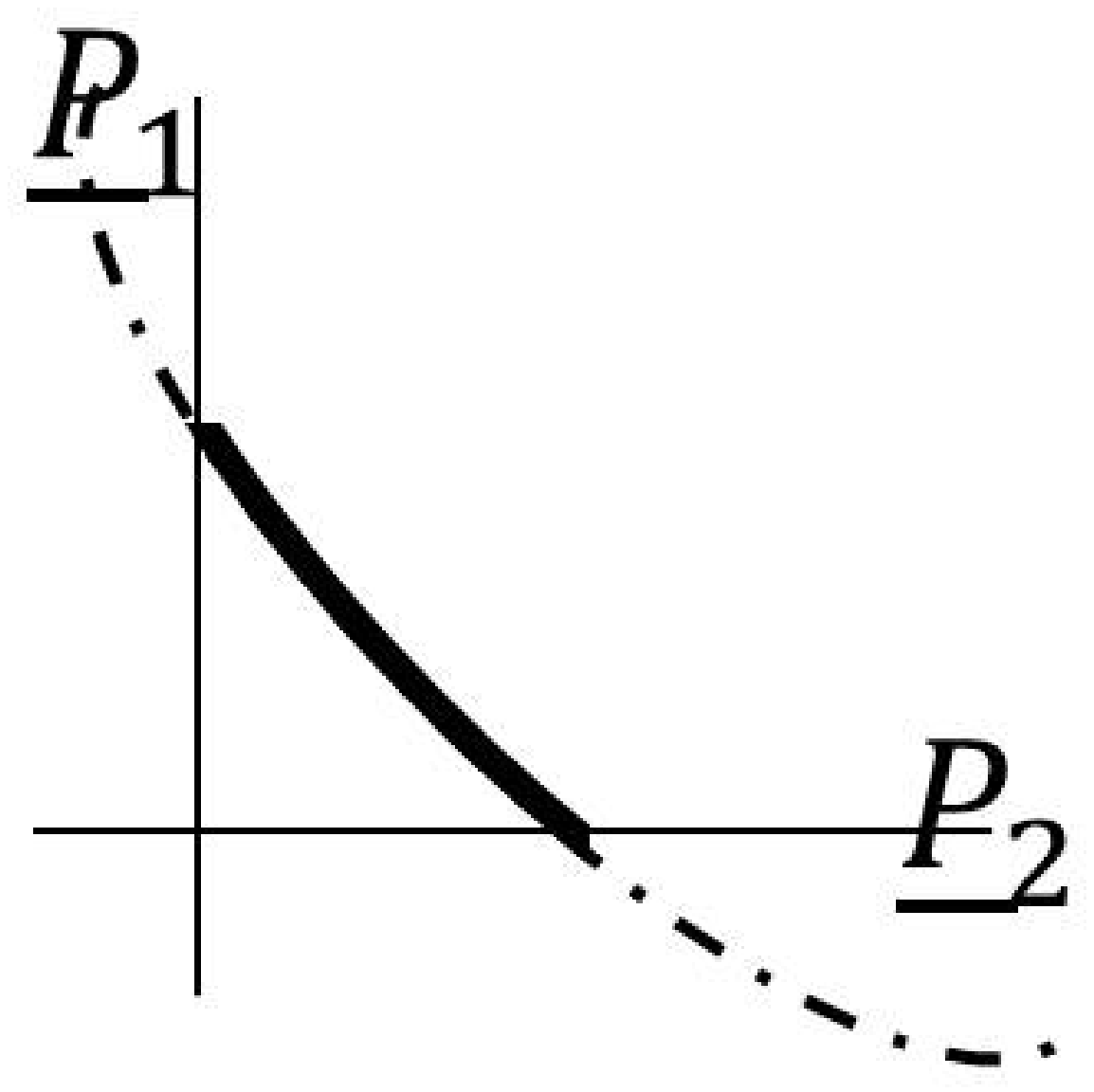}}\hspace{1cm}
\subfigure[Infeasible]{ \label{fig:2bus_theta2}
\includegraphics[width=2cm]{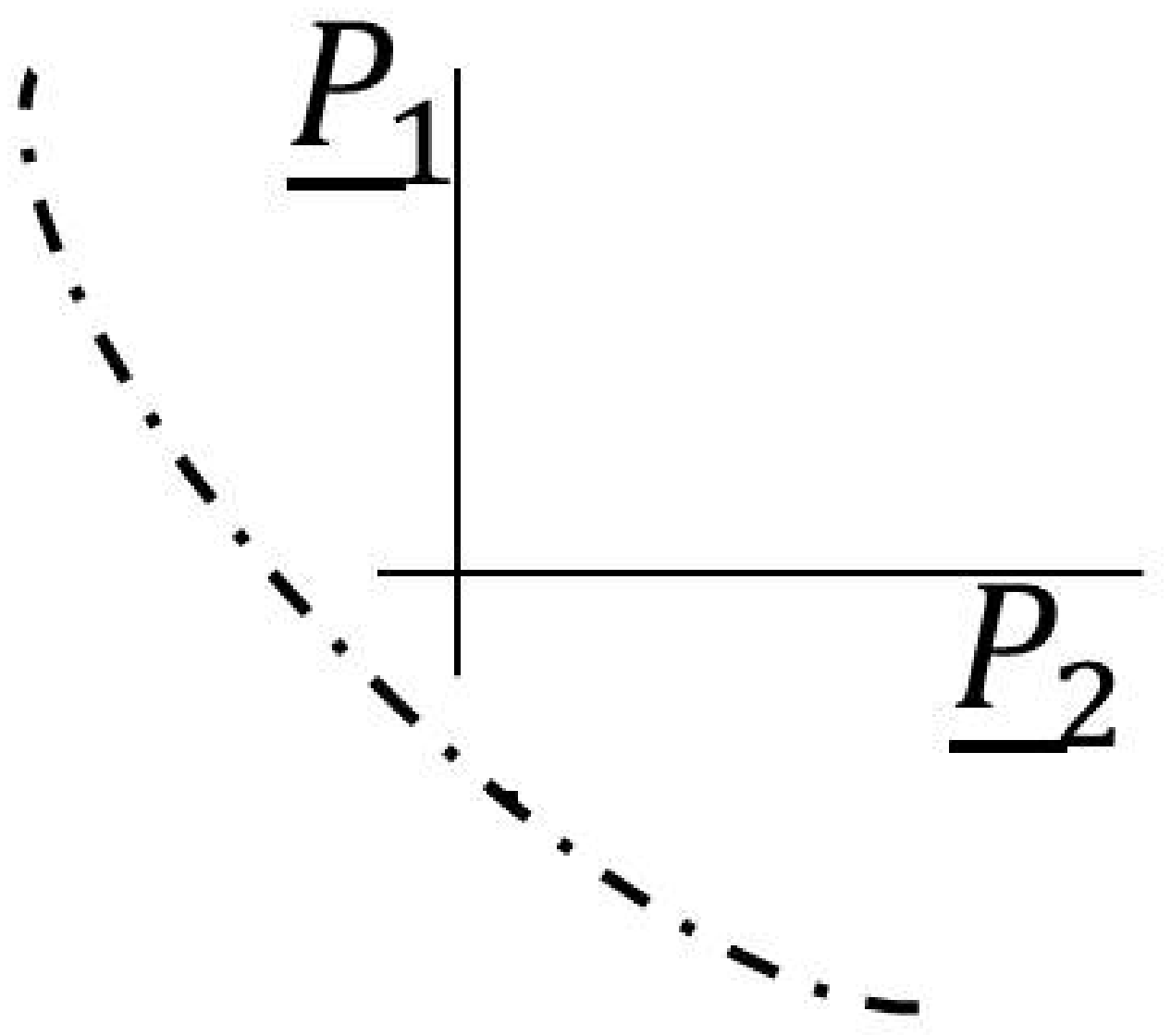}}
\caption{Either $\mc{P} = \mc{O}(\mc{P})=\mathcal{O}(\conv(\mathcal{P}))$ or the injection region is empty.}%
\label{fig:2bus_theta}%
\vspace{-4mm}
\end{figure}

Assume that the angle constraints are such that $\mc{P}_\te=\mc{O}(\conv(\mc{P}_\te))$, implying that every point in $\mc{P}_\te$ is Pareto-optimal. Now, there are two possible scenarios for the injection region $\mc{P}$ as shown in Figure~\ref{fig:2bus_theta}.  In Figure \ref{fig:2bus_theta1}, some of the points of the region $\mc{P}_\te$ remain in $\mathcal P$ and they form the Pareto-front of both $\mc{P}$ and $\conv(\mc{P})$. In Figure \ref{fig:2bus_theta2}, $\mc{P}=\emptyset$ so then $\conv(\mc{P})=\emptyset$ as well. We observe in both cases that $\mc{O}(\mc{P}) =\mathcal{O}(\conv(\mathcal{P}))$. Therefore, we have $\mc{O}(\mc{P})=\mc{O}(\conv(\mc{P}))$ if $\mc{P}_\te = \mc{O}(\conv(\mc{P}_\te))$.

In terms of the line parameters $b_{12}$ and $g_{12}$, the condition $\mc{P}_\te = \mc{O}(\conv(\mc{P}_\te))$ can be written as:
\begin{equation}
\label{eq_n_r_1}
-\tan^{-1} (\frac{b_{12}}{g_{12}}) < \ul{\te}_{12} \leq \ov{\te}_{12} < \tan^{-1} (\frac{b_{12}}{g_{12}}).
\end{equation}
Observe that $\tan^{-1} (b_{12}/g_{12})$ is equal to $45.0^{\circ}$, $63.4^{\circ}$ and $78.6^{\circ}$ for $\frac{b_{12}}{g_{12}}$ equal to $1$, $2$ and $5$, respectively. These numbers suggest that the above condition is very practical. For example, if the inductance of a transmission line is larger than its resistance, the above requirement is met if $| \ul{\te}_{12}|, |\ov{\te}_{12}|< 45.0^{\circ}$. It is noteworthy that an assumption
\begin{equation}
\label{eq_n_r_11}
-\tan^{-1} (\frac{g_{12}}{b_{12}}) < \ul{\te}_{12} \leq \ov{\te}_{12} < \tan^{-1} (\frac{g_{12}}{b_{12}})
\end{equation}
is made in Chapter 15 of \cite{Baldick06}, under which a practical optimization can be convexified (after approximating the power balance equations). However, our condition (\ref{eq_n_r_1}) is less restrictive  than (\ref{eq_n_r_11}) if $b>g$. To understand the reason, note that the value $\frac{g_{12}}{b_{12}}$ is around $0.1$  for a typical transmission line at the transmission level of the network \cite{Baldick06}  . Now, our condition allows for an angle difference as high as $80^o$ while the condition reported in \cite{Baldick06} confines the angle to $6^o$.

\subsection{General Tree Networks} \label{sec:O}
In this section, we study general tree networks with local angle constraints and global bus power constraints. For every bus $i\in\mathcal V$, let $\ov{V_i}$ denote the fixed voltage magnitude $|V_i|$. Given an edge $(i,k)\in\mathcal E$,
 assume that the angle difference $\theta_{ik}$ belongs to the interval $[\ul{\te}_{ik},\ov{\te}_{ik})]$, where $\ul{\te}_{ik}\in[-\pi,0]$ and $\ov{\te}_{ki}\in[0,\pi]$. Define the angle-constrained flow region for the line $(i,k)$ as
\begin{equation}
\nonumber
\mc{F_\te}_{ik}=\{(P_{ik},P_{ki}): \ul{\te}_{ik} \leq \te_{ik} \leq \ov{\te}_{ik},|V_i|=\ov{V}_i, |V_k|=\ov{V}_k\}
\end{equation}
The angle-constrained injection region can be expressed as {$\mc{P}_\te= \bd{A} \mc{F}_\te$}, where $\mc{F}_\te=\prod_{(i,k)\in \mc{E}} \mc{F_\te}_{ik}$. Following the insight from the last subsection, we make the following practical assumption
\begin{equation}
\label{eq_n_r_e1}
-\tan^{-1} (\frac{b_{ik}}{g_{ik}}) < \ul{\te}_{ik} \leq \ov{\te}_{ik} < \tan^{-1} (\frac{b_{ik}}{g_{ik}}),\quad \forall (i,k)\in\mathcal E.
\end{equation}
This ensures that  all points in the flow region of every line are  Pareto optimal. As will be shown later, this assumption leads to the invertibility of the mapping from the injection region $\mc{P}_\te$ to the flow region $\mc{F}_\te$, or equivalently the uniqueness of the solution of every power flow problem.

Assume that the power injection $P_i$  must be within the interval $[\ul{P}_i, \ov{P}_i]$ for every $i\in\mathcal V$. To account for these constraints, define the hyper-rectangle $\mc{P}_P=\{\bd{p}: \ul{\p} \leq  \p \leq \ov{\p} \}$, where $\ul{\p}=(\ul{P}_1,\dots,\ul{P}_n)$ and $\ov{\p}=(\ov{P}_1,\dots,\ov{P}_n)$. The injection region $\mathcal P$ is then equal to $\mc{P}_\te \cap \mc{P}_P$. In what follows, we present the main result of this section.

\begin{thm} \label{lem:angle}
Suppose that $\mathcal P$ is a non-empty set. Under the assumption (\ref{eq_n_r_e1}), the following statements hold:
\begin{enumerate}
\item For every injection vector $\bd{p} \in \mc{P}$, there exists a unique flow vector $\bd{f} \in \mc{F}$ such that $\bd{A} \bd{f} =\bd{p}$.
	\item $\mathcal{P}=\mathcal{O}(\mathcal{P})$.
 \item   $\mathcal{O}(\mathcal{P})=\mathcal{O}(\conv(\mathcal{P}))$.
\end{enumerate}
\end{thm}

In order to prove this theorem, the next lemma is needed.
\begin{lem} \label{lem:angle1}
Under the assumptions of Theorem~\ref{lem:angle}, the relation $\mathcal{O}(\mathcal{P})=\mathcal{O}(\conv(\mc{P}_\te) \cap \mc{P}_P)$ holds.
\end{lem}

The proof of this lemma is provided in the appendix.  Using this lemma, we  prove Theorem~\ref{lem:angle} in the sequel.

{\it Proof of Part 1:}  Given $\bd{p} \in \mc{P}$, consider an arbitrary leaf vertex  $k$. Assume that $i$ is the parent of bus $k$. Since $k$ is a leaf, we have $P_{ki}=P_k$, and  subsequently $P_{ik}$ can be uniquely determined using the relation  $\mathcal{F}_{jk}= \mathcal{O}(\mathcal{F}_{jk})$. One can continue this procedure for every leaf vertex and then go up the tree to determine the flow along each line in every direction.

{\it Proof of Part 2:} Since $\mc{P}$ is a subset of $\mc{P}_\te$, it is enough to show that $\mathcal{P_\te}=\mathcal{O}(\mathcal{P_\te})$. To prove this, the first observation is that  $\mc{F}_\te=\mc{O}(\mc{F}_\te)=\mc{O}(\conv(\mc{F}_\te))$. Given a point $\bd{p} \in \mc{P}_\te$, let $\bd{f} \in \mc{F}_\te$ be the unique flow vector such that $\bd{A} \bd{f}= \bd{p}$. There exist strictly positive numbers $\{c_{ik}, (i,k) \in \mc{E}\}$ such that $\bd{f}$ is the optimal solution to the following optimization problem
\begin{equation}
\label{eqn:flow_opt_1}
\bd{f} = \arg \min_{\tilde{\bd{f}} \in \mc{F}_\te}  \sum_{(i,k) \in \mc{E}} c_{ik} \tilde{P}_{ik} .
\end{equation}
Since minimizing a strictly increasing function gives rise to a Pareto point, it is enough to show that there exists a set of  positive constants $c_1,c_2,\dots,c_n$ such that the optimal solution of the above optimization does not change if its objective function (\ref{eqn:flow_opt_1}) is replaced by $\sum_{(i,k) \in \mc{E}} c_i \tilde{P}_{ik}= \sum_{i=1}^n c_i \tilde P_i$. Since $\mc{F}_\te$ is a product space, we can multiply any pair $(c_{ik}, c_{ki})$ by a positive constant, and $\bd{f}$  still remains  an optimal solution. Assume that the tree is rooted at $1$. Let $i$ be a leaf of the tree and consider the path from $1$ to $i$. Without loss of generality, assume that the nodes on the path are labeled as $1,2,\dots,i$. By setting $c_1$ as $c_{12}$, one can define $c_2,\dots,c_i$ according to the following recursion
\begin{equation}
\nonumber
c_k=c_{k-1} \frac{c_{k,k-1}}{c_{k-1,k}},
\end{equation}
where $k$ ranges from $2$ to $i$. After defining $c_1,\dots,c_i$, we remove all lines of the path 1--$i$ from the network. This creates $i$ disconnected subtrees of the network rooted at $1,\dots,i$. For each of the subtrees with more than $1$ node, one can repeat the above cost assignment procedure until  $c_1,\dots,c_n$ have all been constructed. This completes the proof.

{\it Proof of Part 3:} For notational simplicity, denote $\conv(\mc{P}_\te) \cap \mc{P}_P$ as $\mc{S}$. To prove this part, we use the relation
\begin{equation}
\label{eq_r_eq1}
\mc{P}\subseteq \conv(\mc{P})\subseteq \mc{S}
\end{equation}
and the result of Lemma~\ref{lem:angle1}, i.e.,
\begin{equation}
\label{eq_r_eq2}
\mc{O}(\mc{P})=\mc{O}( \mc{S})
\end{equation}
The first goal is to show the relation $\mc{O}(\mc{P}) \subseteq \mc{O}(\conv(\mc{P}))$ by contradiction. Consider a vector $\p \in \mc{O}(\mc{P})$ such that $\p \notin \mc{O}(\conv(\mc{P}))$. There exists a vector $\p' \in \mc{O}(\conv(\mc{P}))$ such that $\p' \leq \p$ with strict inequality in at least one coordinate. Hence, it follows from (\ref{eq_r_eq1}) that $\p$ is not a Pareto point of $\mc{S}$, while it is a Pareto point of $\mathcal{P}$. This contradicts (\ref{eq_r_eq2}). To prove the converse statement $\mc{O}(\conv(\mc{P}))  \subseteq \mc{O}(\mc{P})$, consider a point $\p \in \mc{O}(\conv(\mc{P}))$. In light of (\ref{eq_r_eq1}), $\p$ belongs to $\mc{S}$. If $\p \in \mc{O}(\mc{S})$, then $\p \in \mc{O}(\mc{P})$ due to (\ref{eq_r_eq2}). If $\p \notin  \mc{O}(\mc{S})$, then there must exist a point $\p' \in \mc{O}(\mc{S})=\mc{O}(\mc{P})$ such that $\p' \leq \p$ with strict inequality in at least one coordinate. This implies  $\p' \in \mc{P}$ and consequently $\p' \in \conv(\mc{P})$, which contradicts $\p \in \mc{O}(\conv(\mc{P}))$.~\hfill$\blacksquare$

{

\subsection{Numerical Algorithms for Convexification}\label{sec:app_fixed}
 The goal of this part is to understand how the results of the preceding subsection can be used to numerically solve an optimization with the feasible set $\mathcal P$. 
The relation $\mathcal{O}(\mathcal{P})=\mathcal{O}(\conv(\mathcal{P}))$ derived before states that the minimization of an increasing function over either the nonconvex set $\mathcal{P}$ or the convexified counterpart $\conv(\mathcal{P})$ leads to the same solution.}{ However, employing a numerical algorithm to minimize a function directly over $\conv(\mathcal{P})$ is difficult due to the lack of efficient algebraic representations of $\conv(\mathcal{P})$.  

To address this issue, one can decompose $\mathcal P$ as $\mc{P}_\te \cap \mc{P}_P$ and then use the fact that  $\conv(\mc{P}_\te)$ and $\conv(\mc{P}_P)$ both have simple algebraic representations. Consider the OPF problem
\begin{subequations}
\label{eq_nr_11}
\begin{align}
\min\quad & \sum_{i\in\mathcal V} f_i(P_i) \\
\mbox{subject to } \quad &  \ul{P}_i \leq P_i \leq \ov{P}_i,\quad \hspace{0.7cm}i\in\mathcal V\\
& P_i=\sum_{k\sim i} P_{ik},\qquad i\in\mathcal V\\
& (P_{ik},P_{ki})\in \mathcal F_{\theta_{ik}},\quad  \hspace{1cm}(i,k)\in\mathcal E; \label{eqn:eq_nr_11_last}
\end{align}
\end{subequations}
where  $f_i(\cdot)$ is both monotonically increasing and convex for every $i\in\mathcal V$, $\ul{P}_i$ and $\overline{P}_i$ are the bus active power lower and upper bounds.  Algebraically, the set $\mathcal{F}_{\theta_{ik}}$ can be represented as  all vectors $(P_{ik},P_{ki})$ satisfying the relations
\begin{equation}
\label{eq_ne_1}
\bigg\|\left[\begin{array}{cc} b_{ik}&-g_{ik}\\ -b_{ik}& -g_{ik}\end{array}\right]^{-1}\left[\begin{array}{cc}P_{ik}-|V_i|^2g_{ik}\\P_{ki}-|V_k|^2g_{ik}\end{array}\right]\bigg\|_2=|V_i||V_k|
\end{equation}
and
\begin{equation}
\nonumber
P_{ki}\leq \ul{P}_{ki}+\frac{\ov{P}_{ki}-\ul{P}_{ki}}{\ov{P}_{ik}-\ul{P}_{ik}} ({P}_{ik}-\ul{P}_{ik}), 
\end{equation}
 $(\ul{P}_{ik},\ul{P}_{ki})$ and $(\ov{P}_{ik},\ov{P}_{ki})$ denote the flows associated with line $(i,k)\in\mathcal E$ for $\theta_{ik}$ equal to $\ul{\theta}_{ik}$ and $\ov{\theta}_{ik}$, respectively.
 The set $\mc{P}_\theta$ is a linear mapping from the product flow region $\mc{F}_\theta$, and hence also has an algebraic representation. The convex hull of $\mathcal{F}_{\theta_{ik}}$ (similarly $\mc{F}_\theta$, $\mc{P}_\theta$) is obtained by changing the equal to sign in \eqref{eq_ne_1} to a less than or equal to sign. If we replace \eqref{eqn:eq_nr_11_last} with $(P_{ik},P_{ki})\in \conv(\mathcal F_{\theta_{ik}})$, a convex problem is obtained. This {\it convexified OPF} can be solved as an second order cone problem (SOCP) efficiently in polynomial time \cite{Lavaei11c,opt_j_4,Gan12}. Since any SOCP problem can be written as a semi-definite programming (SDP) problem, the convexified problem  can also be interpreted as an SDP \cite{Zhang12,Lam12b}.
 
 The feasible region of the convexified problem is not $\conv(\mathcal{P})=\conv(\mc{P}_\te \cap \mc{P}_P)$ but instead $\conv(\mc{P}_\te) \cap \mc{P}_P$. In general, the convex hull operation and the intersection operation do not commute. However, by Lemma~\ref{lem:angle1}, the Pareto fronts of $\conv(\mathcal{P})$ and $\conv(\mc{P}_\te) \cap \mc{P}_P$ are identical (see Fig. \ref{fig:conv_inter}). Therefore if the original problem in \eqref{eq_nr_11} is feasible, optimizing over $\conv(\mathcal{P})$ or $\conv(\mc{P}_\te) \cap \mc{P}_P$ yield the same solution. Moreover, 
 \begin{itemize}
   \item If the solution of the convexified OPF is not a feasible point of OPF, then the original OPF problem is infeasible.
   \item If the solution of the convexified OPF is  a feasible point of OPF, then it is a globally optimal solution of the original OPF problem as well.
 \end{itemize}     
 \begin{figure}
 \centering
 \subfigure[$\conv(\mathcal{P})$]{
 \includegraphics[scale=0.2]{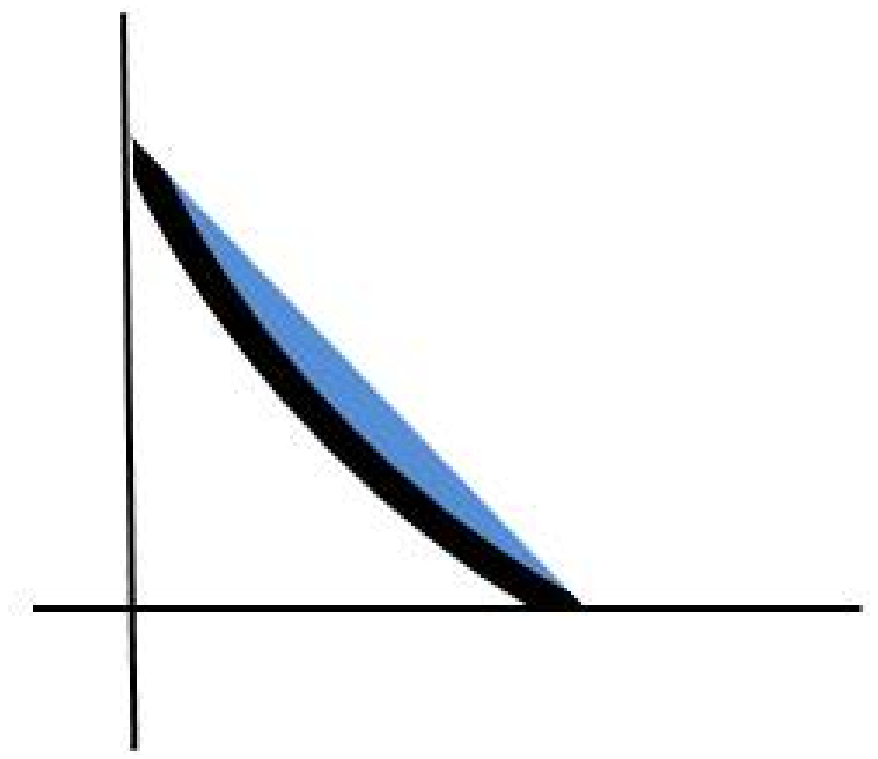}}
 \subfigure[$\conv{\mc{P}_\theta}\cap \mc{P}_P$]{
 \includegraphics[scale=0.2]{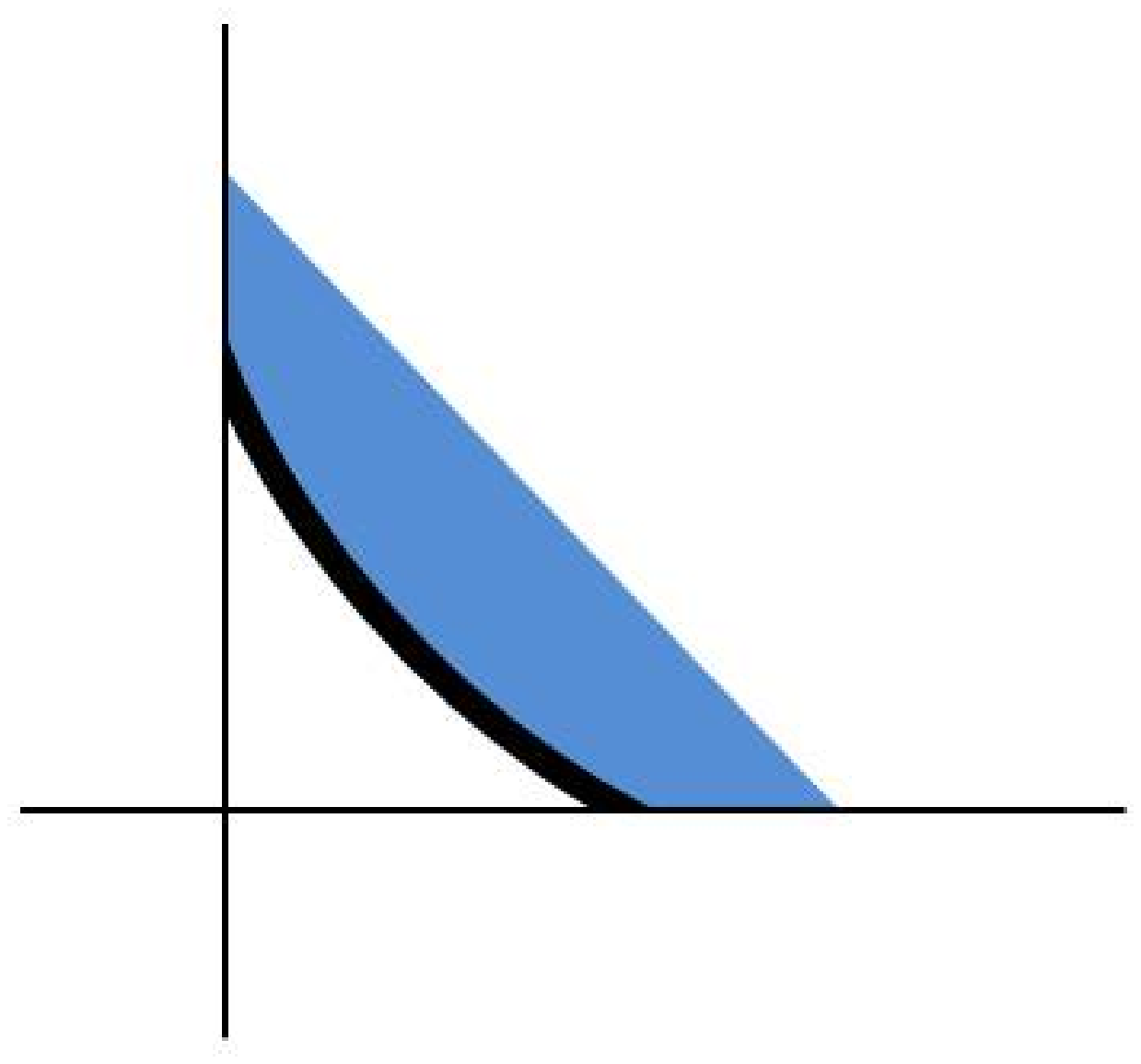}}
 \caption{The sets $\conv(\mathcal{P})$ and $\conv{\mc{P}_\theta}\cap \mc{P}_P$ for a two bus network. The sets are different, but they share the same Pareto Front. }
 \label{fig:conv_inter}
 \end{figure}
}

\subsection{Nonnegative Locational Marginal Prices}\label{sec:lmp}
The economical interpretation of the locational marginal price at a given bus is the increase in the optimal generation cost of adding one unit of load to that particular bus. We will show that this corresponds to the usual practice of defining the LMPs as the Lagrangian multiplier of the power balance equation. The main objective of this part is to prove that LMPs are always nonnegative for a tree network under assumption~(\ref{eq_n_r_e1}).

 We formalize the definition of LMPs for a nonconvex OPF problem in the sequel. Our definition will be consistent with the existing ones for a linearized OPF (such as DC OPF) \cite{kirschen04}. Given a vector $\boldsymbol\epsilon=\left[\begin{array}{ccccc} \epsilon_1&\cdots&\epsilon_n\end{array}\right]$, let $f^*_{\boldsymbol\epsilon}$ denote the optimal objective value of the OPF problem~(\ref{eq_nr_11}) after perturbing its load vector $\bold p_D=\left[\begin{array}{ccccc} P_{D_1}&\cdots&P_{D_n}\end{array}\right]$ as $\bold p_D+\boldsymbol\epsilon$. To have a meaningful definition of LMPs, we must assume that $f^*_{\boldsymbol\epsilon}$ is differentiable at the point $\boldsymbol\epsilon=0$.   A set  $
\{\lambda_1,...,\lambda_n\}$ is called the set of LMPs for buses $1,2,..,n$ if $\sum_{i\in\mathcal V} \lambda_i \epsilon_i$ is the first-order approximation of $
f^*_{\boldsymbol\epsilon}-f^*_{0}
$. As can be seen from this definition, the set of LMPs is unqiue and well-defined.

\begin{thm} \label{lem:price}
Under the  angle assumption~(\ref{eq_n_r_e1}), the following statements hold for every  $i\in\mathcal V$:
\begin{enumerate}
  \item The LMP $\lambda_i$ is equal to the Lagrange multiplier for the power balance equation $P_{G_i}-P_{D_i}=\sum_{k\sim i} P_{ik}$ in the convexified OPF problem.
  \item The LMP $\lambda_i$ is nonnegative.
 \end{enumerate}
\end{thm}

\emph{Proof of Part 1:} Since the angle constraint (\ref{eq_n_r_e1}) does not depend on $\bold p_D$, it follows from the argument made in Section~\ref{sec:app_fixed} that $f^*_{\boldsymbol\epsilon}$  is equal to the optimal value of the convexified OPF problem after perturbing $\bold p_D$ as $\bold p_D+\boldsymbol\epsilon$. Now,  $\lambda_i$ being  the Lagrange multiplier for the power balance equation $P_{G_i}-P_{D_i}=\sum_{k\sim i} P_{ik}$ is an immediate consequence of the well-known sensitivity analysis in convex optimization \cite{boyd04}.

\emph{Proof of Part 2:}  Denote the Lagrange multipliers for the inequalities $ \ul{P}_{G_i} \leq P_{G_i}$ and  $P_{G_i}\leq \ov{P}_{G_i}$ in the convexified problem as $\ul{\lambda}_i$ and $\ov{\lambda}_i$, respectively, for every $i\in\mathcal V$.  We use the superscript "*" to denote the parameters of the (convexified) OPF problem at optimality. It is straightforward to verify that $\lambda_i=f'(P_{G_i}^*)-\ul{\lambda}_i+\ov{\lambda}_i$. To prove the theorem by contradiction, assume that some LMPs are strictly negative. In line with the proof of Lemma~\ref{lem:angle1}, it can be shown that (see (\ref{eq_r_eq4})):
\begin{equation}
\label{eq_nr_12}
\begin{aligned}
(P_{ik}^*,P_{ki}^*)&= \mathop{\arg\min}_{(P_{ik},P_{ki})\in\conv( \mc{F}_{\te_{ik}})} \lambda_i P_{ik}+ \lambda_k P_{ki}\\
&=\mathop{\arg\min}_{(P_{ik},P_{ki})\in  \mc{F}_{\te_{ik}}}\lambda_i P_{ik}+ \lambda_k P_{ki}
\end{aligned}
\end{equation}
for every $(i,k)\in\mathcal E$. Define $\mathcal T$ as a connected, induced subtree of the network with the maximum number of vertices such that $\lambda_i<0$ for every $i\in\mathcal T$. A node $k\in\mathcal V\backslash \mathcal T$ is called a neighbor of $\mathcal T$ if $(i,k)\in\mathcal E$ for some $i\in\mathcal T$. Define ${\mathcal T}_e$ as the subgraph induced by $\mathcal T$ and its neighbors. For every line $(i,k)\in\mathcal E$, let $\ov{P}_{ik}$ denote the maximum possible flow from node $i$ to node $k$ on this line, i.e., $\ov{P}_{ik}=\max\{P_{ik}| (P_{ik},P_{ki})\in\mc{F_\te}_{ik}\}$.
For every $(i,k)\in\mathcal T_e$, it can be deduced from (\ref{eq_nr_12}) and  the geometric properties of $\conv( \mc{F_\te}_{ik})$ and $\mc{F_\te}_{ik}$ (as studied earlier) that
\begin{equation}
\label{eq_nr_14}
\begin{aligned}
P_{ik}^*=\ov{P}_{ik}\qquad\qquad&\quad\qquad \text{if } i\in\mathcal T, k\not\in\mathcal T\\
P_{ik}^*=\ov{P}_{ik}\quad \text{or}\quad P_{ki}^*=\ov{P}_{ki}&\qquad\quad \text{if } i,k\in\mathcal T
\end{aligned}
\end{equation}
We orient the edges of $\mathcal T_e$ to obtain a direct graph $\vec{\mathcal T}_e$ using the following procedure: for every line $(i,k)\in\mathcal T_e$, the orientation of this edge in $\vec{\mathcal T}_e$ is from vertex $i$ to vertex $k$ if $P_{ik}^*=\ov{P}_{ik}$. Since $\vec{\mathcal T}_e$ is a directed tree, it must have a node $i$ whose in-degree is zero. Due to (\ref{eq_nr_14}), node $i$ must belong to $\mathcal T$. Now, one can write:
\begin{equation}
\label{eq_nr_13}
P_{ik}^*=\ov{P}_{ik},\qquad \forall\ k\sim i
\end{equation}
On the other hand, since $i\in\mathcal T$, we have $\lambda_i<0$ and therefore $P_{G_i}^*=\ul{P}_{G_i}$ (note that $\lambda_i=f'(P_{G_i}^*)-\ul{\lambda}_i+\ov{\lambda}_i$).  This implies that
 \begin{equation}  
\nonumber
\ul{P}_{G_i}=\sum_{k\sim i}\ov{P}_{ik}
\end{equation}
As a result of this equality, if the load $P_{D_i}$ is perturbed as $P_{D_i}+\epsilon$ for a very small number $\epsilon<0$, the OPF problem becomes infinite (because when  generator $i$ operates at its minimum output with $\epsilon=0$,  the flows on all lines connected to bus $i$ will hit their maximum values). This contradicts the differentiability of $f^*_{\boldsymbol\epsilon}$ at $\boldsymbol\epsilon=0$. Hence, $\lambda_i$ must be nonnegative.~\hfill$\blacksquare$

The convexification method proposed in \cite{javad_zero1,Lavaei11c,Bose11} for the OPF problem relies on a load over-satisfaction assumption. In the language of our work, this is tantamount to the non-negativity of the LMPs. As a by-product of Theorem~\ref{lem:price}, we have shown that this load over-satisfaction assumption made in previous papers is satisfied as long as  a practical angle assumption holds.

\section{Variable Voltage Pareto Optimal Points} \label{sec:variable}
So far, we have assumed that all complex voltages in the network have fixed magnitudes. In this section, the results derived earlier will be extended to the case with variable voltage magnitudes under the assumption $\underline P_i=-\infty$ for every  $i\in\mathcal V$. The goal is to study the injection region after imposing the constraints 
\begin{subequations}
\begin{align}
&P_i\leq \ov{P}_i,\qquad \qquad  \quad i\in\mathcal V\\
\label{eq_nr_4}
&\theta_{ik}\in[\ul{\theta}_{ik},\ov{\theta}_{ik}],\qquad    (i,k)\in\mathcal E
\end{align}
\end{subequations}
where $|\ul{\theta}_{ik}|,|\ov{\theta}_{ik}|<90^{\circ}$. Note that the results to be developed next are valid even with explicit line flow constraints.

Given a bus $i\in\mathcal V$, let $\ul{V}_i$ and $\ov{V}_i$ denote the given lower and upper bounds on $|V_i|$.  In vector notation, define $\ul{\bd{v}}=(\ul{V}_1,\dots,\ul{V}_n)$ and $\ov{\bd{v}}=(\ov{V}_1,\dots,\ov{V}_n)$. Given a vector $\tilde{\bold v}\in\R^+$,  define $\mc{P}_{\theta}(\tl{\bd{v}})$ as the angle-constrained injection region  in the case when  the voltage magnitudes are fixed according to $\tl{\bd{v}}$, i.e. $|V_i|=\tl V_i$ for $i=1,...,n$. Let $\mathcal P$ and ${\mathcal P}_{\theta}$ denote the regions for the case with  variable voltage magnitudes. One can write $\mathcal P={\mathcal P}_{\theta}\cap {\mathcal P}_P$, where
\begin{equation}
\label{eq_nr_2}
\mc{P}_{\theta}= \bigcup_{\ul{\bd{v}} \leq \tl{\bd{v}} \leq \ov{\bd{v}}} \mc{P}_{\theta}(\tl{\bd{v}})
\end{equation}
The problem of interest is to compute the convex hull of $\mc{P}_{\theta}$. However, the challenge is that 
 the union operator does not commute with the convex hull operator in general (because the union of two convex sets may not be convex).  In what follows, this issue will be addressed by exploiting the flow decomposition technique introduced in \cite{Lavaei11c}. Let $\mc{H}_2^+$ denote the convex set of $2\times 2$ positive semidefinite Hermitian matrices and $\mathcal H_n$ denote the set of all $n\times n$ Hermitian matrices. Given a matrix $\bold W\in\mathcal H_n$ together with an edge $(i,k)\in\mathcal E$, define:
 \begin{itemize}
\item $W_{ik}$: $(i,k)$ entry of $\bold W$.
\item $\bold W_{ik}$: The $2\times 2$ submatrix of $\bold W$ corresponding to the entries  $(i,i),(i,k),(k,i),(k,k)$. The matrix $\bold W_{ik}$ is called an edge submatrix of $\bold W$.
\end{itemize}
Define also
\begin{equation}
\nonumber
\begin{aligned}
{\mc{H}}_{ik}(\tl{\bd{v}})=\bigg\{ \bd{W} \in \mc{H}:\ &\bold W_{ik}  \in \mc{H}_2^+, W_{ii}=\tl{V}_i^2, W_{kk}=\tl{V}_k^2, \\
& \tan(\ul{\theta}_{ik}) \times \text{Im}(W_{ik}) \leq \text{Re}(W_{ik}),  \\
&  \text{Re}(W_{ik})\leq \tan(\ov{\theta}_{ik}) \times \text{Im}(W_{ik})   \bigg\}
\end{aligned}
\end{equation}
It can be shown that for every matrix $\bold W\in {\mc{H}}_{ik}(\tl{\bd{v}})$ with the property $\text{Rank}(\bold W_{ik})=1$, there exists an angle $\theta_{ik}\in[\ul{\theta}_{ik},\ov{\theta}_{ik}]$ such that
\begin{equation}
\nonumber
\bold W_{ik}=\left[\begin{array}{ccc} \tl{V}_i^2 & \tl{V}_i\tl{V}_k \measuredangle\theta_{ik}\\ \tl{V}_i\tl{V}_k \measuredangle\theta_{ki} & \tl{V}_k^2\end{array}\right]
\end{equation}
Thus, 
\begin{equation}
\nonumber
\begin{aligned}
\mc{F}_{\theta_{ik}}(\tl{\bd{v}})= \big\{\Real(\diag(\bd{W}_{ik} \bd{Y}_{ik}^H)):  \bd{W} &\in {\mc{H}}_{ik}(\tl{\bd{v}}),\\
& \text{Rank}(\bold W_{ik})=1\big\}
\end{aligned}
\end{equation}
where
$\bd{Y}_{ik}=\bma y_{ik} & - y_{ik} \\ -y_{ik} & y_{ik} \ebma$ (see \cite{Zhang12,Lavaei11c}).
The flow region $\mc{F}_{\theta_{ik}}(\tl{\bd{v}})$ can be naturally convexified by dropping its rank constraint. However, the convexified set may not be identical to $\conv(\mc{F}_{\theta_{ik}}(\tl{\bd{v}}))$. We use the notation $\ov{\conv}(\mc{F}_{\theta_{ik}}(\tl{\bd{v}}))$ for the convexified flow region, which is defined as
\begin{equation}
\label{eq_nr_8}
\begin{aligned}
\ov{\conv}(\mc{F}_{\theta_{ik}}(\tl{\bd{v}}))&= \big\{\Real(\diag(\bd{W}_{ik} \bd{Y}_{ik}^H)):  \bd{W} \in {\mc{H}}_{ik}(\tl{\bd{v}})\big\}
\end{aligned}
\end{equation}
The following sets can also be defined in a natural way:
\begin{equation}
\nonumber
\begin{aligned}
\ov{\conv}(\mc{F}_{\theta}(\tl{\bd{v}}))&=\prod_{(i,k) \in \mc{E}} \ov{\conv}(\mc{F}_{\theta_{ik}}(\tl{\bd{v}})),\\
\ov{\conv}(\mc{P}_{\theta}(\tl{\bd{v}}))&=\bold A \ov{\conv}(\mc{F}_{\theta}(\tl{\bd{v}}))
\end{aligned}
\end{equation}
Note that ${\conv}(\cdot)$ and $\ov{\conv}(\cdot)$ were the same if  the angle constraint (\ref{eq_nr_4}) did not exist.

 \begin{lem} \label{lem:union}
Given a vector $\tl{\bd{v}}$, the following relations hold:
\begin{subequations}
\begin{align}
\label{eq_nr_5}
&\mathcal O(\mathcal P_{\theta}(\tl{\bd{v}}))=\mathcal O(\ov{\conv}(\mc{P}_{\theta}(\tl{\bd{v}}))),\\
\label{eq_nr_6}
&\mathcal O(\mathcal P(\tl{\bd{v}}))=\mathcal O(\ov{\conv}(\mc{P}_{\theta}(\tl{\bd{v}}))\cap \mathcal P_P),\\
\label{eq_nr_7}
&\mathop{\bigcup}_{\ul{\bd{v}}\leq\tl{\bd{v}}\leq\ov{\bd{v}}} \ov{\conv}(\mc{P}_{\theta}(\tl{\bd{v}}))=\text{\emph{Convex set}}.
\end{align}
\end{subequations}

\end{lem}

\emph{Proof:} The proofs provided in Section~\ref{sec:O1} for the case with fixed voltage magnitudes can be easily adapted to prove (\ref{eq_nr_5}) and (\ref{eq_nr_6}). Therefore, we only prove (\ref{eq_nr_7}) here. It follows from (\ref{eq_nr_8}) that 
there exists a linear transformation  from ${\mc{H}}_{ik}(\tl{\bd{v}})$ to $\ov{\conv}(\mc{F}_{\theta_{ik}}(\tl{\bd{v}}))$, which is independent of $\tl{\bd{v}}$. We denote this transformation as $\ov{\conv}(\mc{F}_{\theta_{ik}}(\tl{\bd{v}}))= l_{ik} ( {\mc{H}}_{ik}(\tl{\bd{v}}))$. Define ${\mc{H}}(\tl{\bd{v}})$ as $\bigcap_{(i,k)\in\mathcal E} {\mc{H}}_{ik}(\tl{\bd{v}})$ and the function $l(\cdot)$ as the natural extension of $l_{ik}$. Hence, $\ov{\conv}(\mc{F}_{\theta}(\tl{\bd{v}}))= l({\mc{H}}(\tl{\bd{v}}))$.
One can write:
\begin{equation}
\label{eq_nr_9}
\begin{aligned}
\bigcup_{\ul{\bd{v}} \leq \tl{\bd{v}} \leq \ov{\bd{v}}} \ov{\conv}(\mc{P}_{\theta}(\tl{\bd{v}})) &=  \bigcup_{\ul{\bd{v}} \leq \tl{\bd{v}} \leq \ov{\bd{v}}} \bd{A} \ov{\conv}( \mc{F}_{\theta}(\tl{\bd{v}})) \\
&= \bd{A} \bigcup_{\ul{\bd{v}} \leq \tl{\bd{v}} \leq \ov{\bd{v}}} l ({\mc{H}}(\tl{\bd{v}})) \\
&= \bd{A} l \bigg(\bigcup_{\ul{\bd{v}} \leq \tl{\bd{v}} \leq \ov{\bd{v}}} {\mc{H}}(\tl{\bd{v}})\bigg) 
\end{aligned}
\end{equation}
On the other hand, $\bigcup_{\ul{\bd{v}} \leq \tl{\bd{v}} \leq \ov{\bd{v}}} {\mc{H}}(\tl{\bd{v}})$ is a convex set because it consists of all Hermitian matrices $\bold W$ whose entries satisfy certain linear and convex constraints. Due to the convexity of this set as well as the linearity of $\bd{A}$ and $l$, it can be concluded from (\ref{eq_nr_9}) that $\bigcup_{\ul{\bd{v}} \leq \tl{\bd{v}} \leq \ov{\bd{v}}} \ov{\conv}(\mc{P}_{\theta}(\tl{\bd{v}}))$ is convex.~\hfill$\blacksquare$

As pointed out before Lemma~\ref{lem:union},  ${\conv}(\cdot)$ and $\ov{\conv}(\cdot)$ are equivalent if  the angle constraint (\ref{eq_nr_4}) is ignored. In this case, it follows from (\ref{eq_nr_7}) and the relation 
\begin{equation}
\nonumber
\mc{P} _{\theta}\subseteq\bigcup_{\ul{\bd{v}} \leq \tl{\bd{v}} \leq \ov{\bd{v}}} \conv(\mc{P}_{\theta}(\tl{\bd{v}}))\subseteq \conv(\mc{P}_{\theta})
\end{equation}
that $\conv(\mc{P}_{\theta})=\bigcup_{\ul{\bd{v}} \leq \tl{\bd{v}} \leq \ov{\bd{v}}} \conv(\mc{P}_{\theta}(\tl{\bd{v}}))$. In other words, as long as there is no angle constraint, the convex hull operator commutes with the union operator when it is applied to (\ref{eq_nr_2}). Motivated by this observation, define $\ov{\conv}(\mc{P}_{\theta})$ as the convex set $\bigcup_{\ul{\bd{v}} \leq \tl{\bd{v}} \leq \ov{\bd{v}}} \ov{\conv}(\mc{P}_{\theta}(\tl{\bd{v}}))$. We present the main theorem of this section below.

\begin{thm} \label{thm:EOunion}
For a tree network, $\mc{O}(\mc{P}) = \mc{O}(\conv(\mc{P}))=\mathcal{O}(\ov{\conv}(\mc{P}_\te) \cap \mc{P}_P)$.
\end{thm}

\emph{Proof:} Since 
\begin{equation}
\label{eq_nr_10}
\mc{P}\subseteq\conv(\mc{P})\subseteq \ov{\conv}(\mc{P}_\te) \cap \mc{P}_P,
\end{equation}
it suffices to prove that  $\mc{O}(\mc{P}) = \mathcal{O}(\ov{\conv}(\mc{P}_\te) \cap \mc{P}_P)$  (see part 3 of Theorem~\ref{lem:angle}). First, we show that $\mathcal{O}(\ov{\conv}(\mc{P}_\te) \cap \mc{P}_P)\subseteq\mc{O}(\mc{P})$. Consider a vector $\p$ in $\mathcal{O}(\ov{\conv}(\mc{P}_\te) \cap \mc{P}_P)$. By the definition of $\ov{\conv}(\mc{P}_\te)$, $\p \in \mathcal O(\ov{\conv}(\mc{P}_\te(\tl{\bd{v}}))\cap \mc{P}_P)$ for some $\tl{\bd{v}}$. Hence, by Lemma~\ref{lem:union}, $\p \in \mc{O}(\mc{P}(\tl{\bd{v}}))$ and consequently $\p\in\mathcal P$. Now, it follows from (\ref{eq_nr_10}) and $\p\in\mathcal{O}(\ov{\conv}(\mc{P}_\te) \cap \mc{P}_P)$ that $\p\in\mathcal O(\mathcal P)$. 
The relation $\mc{O}(\mc{P})\subseteq \mathcal{O}(\ov{\conv}(\mc{P}_\te) \cap \mc{P}_P)$ can be proved in line with the proof of Lemma \ref{lem:angle1}.~\hfill$\blacksquare$

{

\subsection{Convexification via SOCP and SDP Relaxations}  }

Theorem~\ref{thm:EOunion} presents two relations $\mc{O}(\mc{P}) = \mc{O}(\conv(\mc{P}))$ and $\mc{O}(\mc{P}) =\mathcal{O}(\ov{\conv}(\mc{P}_\te) \cap \mc{P}_P)$. Although the first relation reveals a convexity property of $\mathcal P$, it cannot be used directly to convexify a hard optimization over $\mathcal P$. Instead, one can deploy the second relation for this purpose.
By generalizing the argument made in Subsection~\ref{sec:app_fixed}, we will spell out in this sequel how to convexify an OPF problem with variable voltage magnitudes. Consider the following OPF problem:
\begin{equation}
\nonumber
\begin{aligned}
\min\quad & \sum_{i\in\mathcal V} f_i(P_{G_i}) \\
\mbox{subject to } \quad &P_{G_i} \leq \ov{P}_{G_i},\quad \hspace{1.85cm}i\in\mathcal V\\
& P_{G_i}-P_{D_i}=\sum_{k\sim i} P_{ik},\quad  \hspace{0.4cm} i\in\mathcal V\\
& (P_{ik},P_{ki})\in \mathcal F_{\theta_{ik}}(\tl{\bd{v}}) ,\quad \hspace{0.5cm}(i,k)\in\mathcal E\\
&\ul{\bd{v}} \leq \tl{\bd{v}} \leq \ov{\bd{v}}
\end{aligned}
\end{equation}
Note that $\tl{\bd{v}}$ is a variable of this optimization accounting for the voltage magnitudes to be optimized. To convexify this optimization, it suffuse to replace the constraint $(P_{ik},P_{ki})\in \mathcal F_{\theta_{ik}}(\tl{\bd{v}}) $ with $(P_{ik},P_{ki})\in \ov{\conv}(\mathcal F_{\theta_{ik}}(\tl{\bd{v}})) $.  { Theorem~\ref{thm:EOunion}  guarantees that OPF and convexified OPF will have the same solution. As shown in Section~\ref{sec:app_fixed} for a special case, the convexified OPF is an SOCP problem.  The statement that OPF and its SOCP relaxations have the same global solution has been proven in \cite{Lavaei11c} using a purely algebraic (rather than geometric) technique. This result is closely related to the prior work \cite{javad_zero1,Zhang12,Bose11}, which shows that OPF and its SDP relaxation have the same solution. This SDP relaxation can be obtained by dropping a single rank constraint $\text{Rank}(\bold W)=1$ as opposed to a series of rank constraints on the edge submatrices of $\bold W$. Note that
\begin{itemize}
  \item As shown in \cite{Lavaei11c}, if the feasible sets of the SDP and SOCP relaxations are projected onto the space for bus injections, they will lead to the same reduced feasible set (this property is true only for tree networks). 
  \item Despite of the fact that these two relaxations convexify the injection region in the same way, it is much easier to solve the SOCP relaxation, from a computational point of view. This is due to the total number of variables involved in the optimization.
\end{itemize}  }

\subsection{Inclusion of Lower Bound on Bus Injection}

Theorem~\ref{thm:EOunion} has been developed under the assumption $\underline P_i=-\infty$ for every  $i\in\mathcal V$. The objective of this part is to remove the above assumption by allowing $\underline P_i$ to be any finite number. The main idea behind this generalization is first to   reduce the variable-voltage-magnitude case to a fixed-voltage-magnitude case and then to deploy Theorem~\ref{lem:angle}. However, Theorem~\ref{lem:angle} is based on two assumptions: (i) non-emptiness of the injection region, (ii) an angle condition for each line. Hence, we need to develop the counterparts of these two assumptions for a general case.

\begin{assum} \label{assum1}
For every vector $\tilde{\bold v}$ satisfying the relation $\ul{\bd{v}} \leq \tl{\bd{v}} \leq \ov{\bd{v}}$, the feasible set $ \mc{P}_{\theta}(\tl{\bd{v}})\cap {\mathcal P}_P$ is nonempty.
\end{assum}

\begin{assum} \label{assum2}
For every line $(i,k)\in\mathcal E$ and vector $\tilde{\bold v}$ satisfying the relation $\ul{\bd{v}} \leq \tl{\bd{v}} \leq \ov{\bd{v}}$, the sets $\mathcal F_{\theta_{ik}}(\tl{\bd{v}})$ and $ \ov{\conv}(\mathcal F_{\theta_{ik}}(\tl{\bd{v}}))$ share the same Pareto front. 
\end{assum}  

\begin{thm}
The relations $\mc{O}(\mc{P}) = \mc{O}(\conv(\mc{P}))=\mathcal{O}(\ov{\conv}(\mc{P}_\te) \cap \mc{P}_P)$ hold for a tree network under Assumptions~\ref{assum1} and \ref{assum2}
\end{thm}

\emph{Proof:} Similar to the proof of Theorem~\ref{thm:EOunion}, it suffices to show that $\mathcal{O}(\ov{\conv}(\mc{P}_\te) \cap \mc{P}_P)\subseteq\mc{O}(\mc{P})$. Consider a vector $\p$ in $\mathcal{O}(\ov{\conv}(\mc{P}_\te) \cap \mc{P}_P)$. By the definition of $\ov{\conv}(\mc{P}_\te)$, $\p \in \mathcal O(\ov{\conv}(\mc{P}_\te(\tl{\bd{v}}))\cap \mc{P}_P)$ for some $\tl{\bd{v}}$. By adopting the proof of Lemma~\ref{lem:angle1} for the fixed-voltage-magnitude case $\tl{\bd{v}}$, it can be concluded that $\p \in \mc{O}(\mc{P}(\tl{\bd{v}}))$ and subsequently $\p\in\mathcal P$ (this requires Assumptions~\ref{assum1} and \ref{assum2}). Now, it follows from (\ref{eq_nr_10}) and the assumption $\p\in\mathcal{O}(\ov{\conv}(\mc{P}_\te) \cap \mc{P}_P)$ that $\p\in\mathcal O(\mathcal P)$.~\hfill$\blacksquare$

\section{Case Study} \label{sec:simu}
 The optimization problem of interest is
\begin{subequations} \label{eqn:opt_sim}
\begin{align}
\mbox{minimize } & \sum_{i=1}^n P_i \\
& \underline{V}_i \leq V_i \leq \overline{V}_i \\
& \underline{P}_i \leq P_i \leq \overline{P}_i \\
& P_{ik}  \leq \overline{P}_{ik} \\
& \underline{Q}_i \leq Q_i \leq \overline{Q}_i \\
& \bd{p}+j\bd{q}= \diag(\bd{v} \bd{v}^H \bd{Y}^H), 
\end{align}
\end{subequations}
where $\sum_{i=1}^n P_i$ is the system loss,  $\underline{V}_i =0.95 p.u.$, $\overline{V}_i=1.05 p.u.$, $\overline{P}_{ik}$ are taken from the transmission line data sheets, $\bd{p}=[P_1 \dots P_n]^T$ and $\bd{q}=[Q_1 \dots Q_n]^T$. The network data are obtained from 34-bus and 123-bus IEEE test feeders \cite{testfeeders}, and the bounds on active and reactive power are determined using two methods:
\begin{enumerate}
\item Every bus has some device to provide active and reactive power (for example, solar inverters \cite{solarmanual}). Let $\tl{P}_i$ and $\tl{Q}_i$ be the reported active and reactive power in the test feeder datasets. In each run, $\overline{P}_i$ is chosen uniformly at random in $[0.8 \tl{P}_i, \tl{P}_i]$ ($\tl{P}_i$'s are negative since buses are withdrawing power), $\underline{P}_i$ is chosen uniformly at random in $[\tl{P}_i,1.2 \tl{P}_i]$; similarly for $\overline{Q}_i$ and $\underline{Q}_i$. For 34-bus and 123-bus networks, 1000 runs are performed for each. 

\item Some buses have devices to provide active and reactive power, but some other buses have fixed active and reactive power requirements. In both 34-bus and 123-bus networks, 20\% of buses are chosen randomly to have variable power bounds for each run, and the bounds are chosen as the same in method (a). 
\end{enumerate}

The main idea of the simulations is to show that the convex relaxation is tight. The optimization problem in \eqref{eqn:opt_sim} can be written as (see \cite{Zhang12})
\begin{subequations} \label{eqn:opt_sim_rank}
\begin{align}
\mbox{minimize } & \sum_{i=1}^n P_i \\
& \underline{V}_i^2 \leq W_{ii} \leq \overline{V}_i^2 \\
& \underline{P}_i \leq P_i \leq \overline{P}_i \\
& P_{ik}  \leq \overline{P}_{ik} \\
& \underline{Q}_i \leq Q_i \leq \overline{Q}_i \\
& \bd{p}+j\bd{q}= \diag(\bd W \bd{Y}^H) \\
& \rank{\bd W}= 1. \label{eqn:rank} 
\end{align}
\end{subequations}
To make \eqref{eqn:opt_sim_rank} convex, we remove the rank constraint \eqref{eqn:rank}, and solve the resulting convex problem using SDPT3 (implemented using Yalmip). 
  The results are summarized in Table \ref{tab:sim_1}. The rank relaxation is tight in all the runs we preformed. 

\begin{table}[ht]
\centering
\begin{tabular}{|c|c|c|}
\hline
& 34-bus & 123-bus \\
\hline
Case a) & 1000 & 1000 \\
\hline
Case b) & 1000 & 1000 \\
\hline
\end{tabular}
\caption{Results of simulation. 1000 runs are performed for each case a) and case b). The entries in the table represents the number of runs where the rank relaxation is tight. We observe that the rank relaxation is tight for all the test runs.}
\label{tab:sim_1}
\end{table} 

\section{Conclusion} \label{sec:con}
This paper is concerned with understanding the geometric properties of the injection region of a tree-shaped power network.
 Since this region is characterized by nonlinear equations, a fundamental resource allocation problem, named {\it optimal power flow} (OPF), becomes nonconvex and hard to solve. The objective of this paper is to show that this highly nonconvex region preserves important properties of a convex set and therefore optimizations over this region can be cast as convex programs. To this end, we have focused on the Pareto front of the injection region, i.e., the set of those points in the injection region that are eligible to be a solution to a typical OPF problem. First, we have studied the case when  the voltage magnitude of every bus is fixed at its nominal value. Although the injection region is still nonconvex, we have shown that the Pareto fronts  of this set and its convex hull are identical under various network constraints as long as a practical angle condition is satisfied. This implies that the injection region can be replaced by its convex hull in the OPF problem without changing the global solution. An implication of this result is that to convexify the OPF problem, its nonlinear constraints can be replaced by simple  linear and norm constraints  and still a global solution of the original problem will be attained. The injection region of a power network with variable voltage magnitudes is also studied.
\appendix
\subsection{\textbf{Thermal Constraints of Distribution Networks}}
\begin{figure}[ht]
\centering
\includegraphics[scale=0.17]{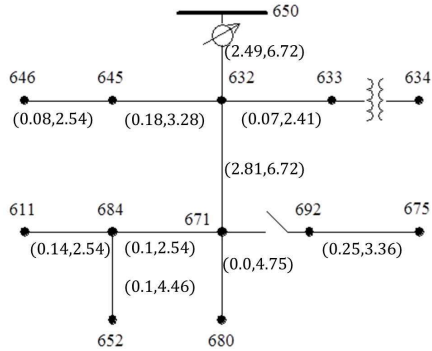}%
\caption{This figure illustrates the angle constraints in a distribution network.}
\label{fig:thermal}%
\vspace{-0.5cm}
\end{figure}

Every transmission line is associated with a current limit, which  restricts the maximum amount of current that can flow through the line. Once this limit and the length of the line are known, we can  convert it first into a $|I|^2 r$ thermal loss constraint and then into an angle constraint on the line. In what follows, we compute some numbers for the 13-bus test feeder system given in \cite{testfeeders}. This system operates at 2.4 KV line to neutral. The angle limits are provided in Figure \ref{fig:thermal}. A pair $(\alpha,\beta)$ is assigned to each line in this figure, where $\alpha$ shows the angle between the two related buses under typical operating conditions (as given in the data) and $\beta$ shows the limit from thermal constraints.

\subsection{Proof of Lemma~\ref{lem:angle1}}

To prove Lemma~\ref{lem:angle1}, we first show that $\mathcal{O}(\mc{S})\subseteq \mathcal{O}(\mc{P})$ (recall that  $\mc{S}=\conv(\mc{P}_\te) \cap \mc{P}_P$). Consider  a point $\x \in \mc{O} (\mc{S})$, and  denote its corresponding line flow from bus $i$ to bus $k$ with $x_{ik}$ for every $(i,k)\in\mathcal E$. Due to the relation  $\mathcal{P}=\mathcal{O}(\mathcal{P})$ derived in Part 2 of Theorem~\ref{lem:angle}, it is enough to prove that $\x\in\mc{P}$. Since $\x$ is a Pareto point of the convex set $\mc{S}$, it is the solution of the following optimization
\begin{equation}
\nonumber
\begin{aligned}
\x = \mathop{\arg\min}_{\p \in \conv(\mc{P_\te})} & \sum_{i=1}^n c_i P_i \\
\mbox{subject to } &  \ul{P}_i \leq P_i \leq \ov{P}_i \,\quad i=1,2,...,n \label{eqn:ul_ov}.
\end{aligned}
\end{equation}
for some positive vector $(c_1,\dots,c_n)$. To simplify the proof, assume that all entries of this vector are strictly positive (the idea to be presented next can be adapted to tackle the case with some zero entries).
By the duality theory, there exist nonnegative Lagrange multipliers $\ul{\la}_1,...,\ul{\la}_n$ and $\ov{\la}_1,...,\ov{\la}_n$ such that
\begin{equation}
\nonumber
\begin{aligned}
\x = \mathop{\arg\min}_{\p \in \conv(\mc{P_\te})} & \sum_{i=1}^n (c_i+\ov{\la}_i-\ul{\la}_i )P_i -\ov{\la}_i \ov{P}_i + \ul{\la}_i \ul{P}_i
\end{aligned}
\end{equation}
or equivalently
\begin{equation}
\label{eq_r_eq3}
\begin{aligned}
\x = \mathop{\arg\min}_{\p \in \conv(\mc{P_\te})} & \sum_{i=1}^n \left(\bar c_i \sum _{k\in\mathcal V:\ k\sim i} P_{ik} \right)
\end{aligned}
\end{equation}
where $\bar c_i=c_i+\ov{\la}_i-\ul{\la}_i $ for every $i\in\mathcal V$. By complementary slackness, whenever $\bar c_i$ is less than or equal to zero,  the multiplier $\ul{\la}_i$ must be strictly positive.  Therefore
\begin{equation}
\label{eq_r_eq5}
 x_i=\ul{P}_i\quad \text{whenever}\quad \bar c_i\leq 0
\end{equation}
 On the other hand, since $\conv(\mc{P_\te})= \bd{A} \conv(\mc{F_\te})$ and  $\mc{F_\te}=\prod_{(i,k) \in \mc{E}} \mc{F}_{\te_{ik}}$, it results from (\ref{eq_r_eq3}) that
\begin{equation}
\label{eq_r_eq4}
(x_{ik},x_{ki})= \mathop{\arg\min}_{(P_{ik},P_{ki})\in\conv( \mc{F_\te}_{ik})} \bar c_i P_{ik}+ \bar c_k P_{ki}
\end{equation}
for every $(i,k)\in\mathcal E$. In order to prove $\x\in\mc{P}$, it suffices to show that $(x_{ik},x_{ki})\in \mc{F_\te}_{ik}$. Notice that if either $\bar c_i>0$ or $\bar c_k>0$, then it can be easily inferred from (\ref{eq_r_eq4}) that $(x_{ik},x_{ki})\in \mc{F_\te}_{ik}$. The challenging part of the proof is to show the validity of this relation in the case when $\bar c_i,\bar c_k\leq 0$. Consider an arbitrary vector $\y$ (not necessarily distinct from $\x$) belonging to $\mc{P}$. Since $(y_{ik},y_{ki})\in \mc{F_\te}_{ik}$, it is enough to prove that   $(x_{ik},x_{ki})=(y_{ik},y_{ki})$ whenever $\bar c_i,\bar c_k\leq 0$. This will be shown below.

Consider an edge $(i,k)\in\mathcal E$ such that $\bar c_i,\bar c_k\leq 0$. There exists at least one connected, induced subtree of the network including the edge $(i,k)$ with the property that $\bar c_r\leq 0$ for every vertex $r$ of this subtree. Among all such subtrees, let $\mathcal G$ denote the one with the maximum number of vertices. We define two types of nodes in $\mc{G}$. A node $r\in\mc{G}$ is called a  boundary node of $\mc{G}$ if either it is connected to some node $l\in\mc{V}\backslash\mc{G}$ or it is a leaf of the tree. We also say that a node $r \in\mc{V}\backslash\mc{G}$ is a neighbor of $\mc{V}$ if it is connected to some node in $\mc{V}$. By (\ref{eq_r_eq5}), if $r$ is a node of $\mc{G}$, then ${y}_r\geq \ul{P}_r=x_r$. Without loss of generality, assume that the tree is rooted at a boundary node of $\mc{G}$, namely node~1.

\begin{figure}[ht]
\centering
\psfrag{Pjk}{$P_{ml}$}
\psfrag{Pkj}{$P_{lm}$}
\psfrag{(xk,xj)}{$(x_{lm},x_{ml})$}
\subfigure[]{\includegraphics[scale=0.4]{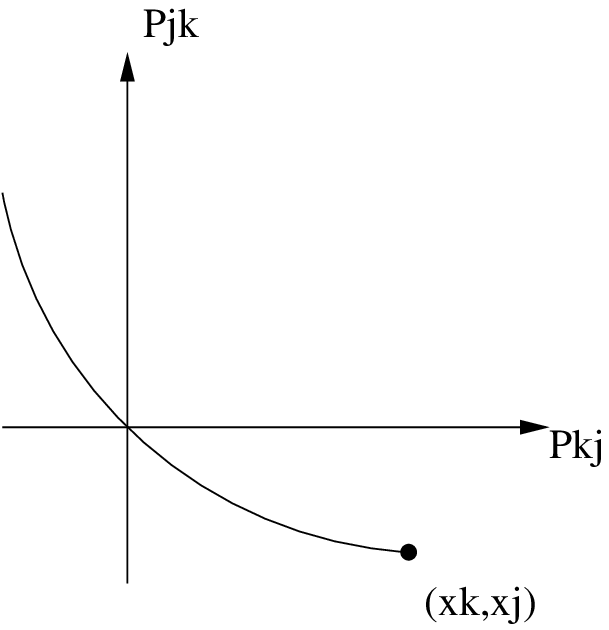}%
\label{fig:Xkj}}
\psfrag{Fjk}{$P_{lr}$}
\psfrag{Fkj}{$P_{rl}$}
\psfrag{(xk,xj)}{$(x_{rl},x_{lr})$}
\psfrag{(pk,pj)}{$(y_{rl},y_{lr})$}\hspace{1.5cm}
\subfigure[]{
\includegraphics[scale=0.4]{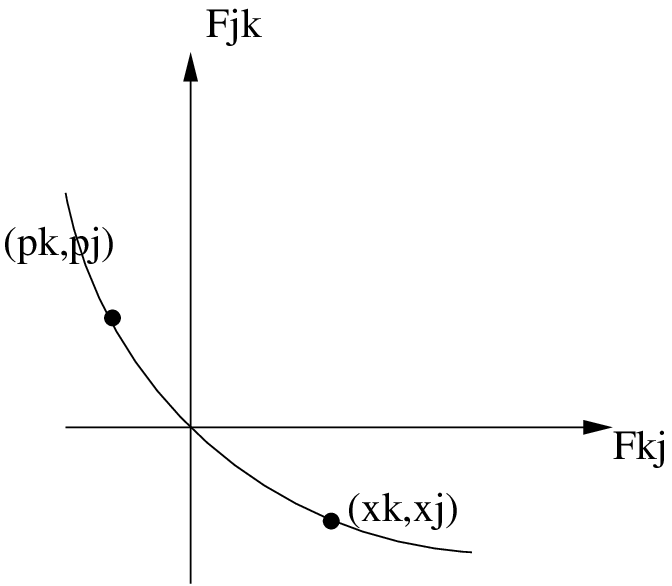}%
\label{fig:PXik}}
\caption{Figure (a) shows the flow region for the line $(l,m)$, where $(x_{lm},x_{ml})$ lies at its lower right corner due to $c_l \leq 0$ and $c_m > 0$. Figure (b) shows the flow region for the line $(r,l)$ to illustrate that $x_{rl} \geq y_{rl}$ (due to $\mc{F}_{\te_{rl}} = \mc{O} (\conv(\mc{F}_{\te_{rl}}))$
and $y_{lr} \geq x_{lr}$).  }
\label{fig:PX2}
\end{figure}

Consider an edge $(r,l)$ of the subtree $\mathcal G$  such that node $l$ is a leaf of  $\mathcal G$ and node $r$ is its parent. First, we want to prove that $y_{lr} \geq x_{lr}$. To this end, consider two possibilities. If $l$ is a leaf of the original tree, then the inequality (\ref{eq_r_eq5}) yields $y_{lr}=y_l \geq \ul{P}_l= x_l=x_{lr}$. As the second case, assume that $l$ is not a leaf of the original tree. Let $m$ denote a neighbor of $\mc{G}$ connected to $l$. By analyzing the flow region for the line $(l,m)$  as depicted in Figure \ref{fig:Xkj},
it follows from (\ref{eq_r_eq4}) and the inequalities $c_l \leq 0, c_m >0$  that $(x_{lm},x_{ml})$ is at the lower right corner of $\mc{F}_{\te_{lm}}$. Thus, $x_{lm} \geq y_{ml}$ because of $(y_{lm},y_{ml})\in \mc{F}_{\te_{lm}}$. Let $\mc{N}_{l}$ denote the set of all nodes connected to $l$ that are neighbors of $\mc{G}$. One can derive the inequality   $x_{lm} \geq y_{lm}$ for every $m \in \mc{N}_{l}$. Combining this set of inequalities with  $x_{l}=\ul{P}_l \leq y_{l}$ or equivalently
\begin{equation}
\nonumber
x_l=\sum_{m \in \mc{N}_l} x_{lm} + x_{lr} \leq \sum_{m \in \mc{N}_l} y_{lm} + y_{lr}=y_l,
\end{equation}
yields that $ y_{lr} \geq x_{lr}$. As illustrated in Figure \ref{fig:PXik}, this implies that    $x_{rl} \geq  y_{rl}$.
This line of argument can be pursued until node~$1$ of the tree is reached. In particular, since node $1$ is assumed to be a boundary node of $\mc{G}$, it can be shown by induction that $x_{1l} \geq y_{1l}$ for every node $l$ such that $(1,l)\in\mathcal E $. On the other hand,
\begin{equation}
\nonumber
\sum_{l\in\mathcal V:\ l\sim1}x_{1l}=x_1=\ul{P}_1\leq y_1=\sum_{l\in\mathcal V:\ l\sim1}y_{1l}
\end{equation}
 Therefore, the equality $x_{1l}=y_{1l}$ must hold for every $l\sim 1$. By propagating this equality down the subtree $\mathcal G$, we obtain that $x_{ik}=y_{ik}$ and $x_{ki}=y_{ki}$. This completes the proof of the relation $\mathcal{O}(\mc{S})\subseteq \mathcal{O}(\mc{P})$.

In order to complete the proof of the lemma, it remains to show that $\mathcal{O}(\mc{P})\subseteq\mathcal{O}(\mc{S})$. To this end, assume by contradiction that there is a point $\p \in \mc{O}(\mc{P})$ such that $\p \notin \mc{O}(\mc{S})$. In light of $\mc{P}\subseteq\mc{S}$, there exists $\p' \in \mc{O}(\mc{S})$ such that $\p' \leq \p$ with strict inequality in at least one coordinate. However,  since $\mathcal{O}(\mc{S})\subseteq \mathcal{O}(\mc{P})$, $\p'$ belongs to $\mc{P}$. This contradicts the assumption $\p \in \mc{O}(\mc{P})$.

\bibliography{mybib}
\bibliographystyle{IEEEtran}
\end{document}

%% file: header.tex
\usepackage{cite}
\usepackage{amsmath, amsthm, amssymb,graphicx,bbm}
\usepackage{verbatim}
\usepackage{psfrag}
\usepackage{algorithm2e}
\usepackage{bm}
\usepackage{color}
\usepackage[tight,footnotesize]{subfigure}

\newcommand{\R}{\mathbb{R}}

\newcommand{\bd}{\mathbf}

\newcommand{\tl}{\tilde}

\newcommand{\x}{\bd{x}}
\newcommand{\y}{\bd{y}}

\newcommand{\ov}{\overline}
\newcommand{\ul}{\underline}
\newcommand{\la}{\lambda}

\newcommand{\C}{\mathbb{C}}

\newcommand{\bma}{\begin{bmatrix}}
\newcommand{\ebma}{\end{bmatrix}}

\newcommand{\nn}{\nonumber}
\newcommand{\te}{\theta}

\newcommand{\mc}{\mathcal}
\newcommand{\p}{\mathbf{p}}

\DeclareMathOperator{\rank}{rank}

\DeclareMathOperator{\diag}{diag}
\DeclareMathOperator{\Real}{Re}

\DeclareMathOperator{\conv}{conv}

\DeclareMathOperator{\conj}{conj}

\newtheorem{thm}{Theorem}
\newtheorem{assum}{Assumption}

\newtheorem{lem}[thm]{Lemma}